\documentclass[11pt]{amsart}

\usepackage[T1]{fontenc}
\usepackage{amssymb,mathrsfs}
\usepackage{palatino,microtype}
\usepackage{enumitem,booktabs,graphicx,float}
\usepackage{tikz}
\usetikzlibrary{arrows.meta,patterns}
\usepackage{pgfplots}
\usepgfplotslibrary{groupplots,fillbetween}
\pgfplotsset{compat=1.18}
\usepackage[colorlinks=true,linkcolor=blue,
            citecolor=blue,urlcolor=blue]{hyperref}
\usepackage{placeins}
\usepackage{xcolor}
\usepackage{mathtools,array}

\DeclareFontFamily{U}{mathx}{}
\DeclareFontShape{U}{mathx}{m}{n}{<-> mathx10}{}
\DeclareSymbolFont{mathx}{U}{mathx}{m}{n}
\DeclareMathAccent{\widecheck}{0}{mathx}{"71}

\hypersetup{pdftitle={Weak Lp integrability of functions with Fourier support in fractal subsets of finite-type manifolds},
  pdfauthor={Thanh-Tinh Cao and Minh-Quy Pham},
  pdfsubject={Fourier support, Hausdorff dimension, Lorentz spaces, finite type}}

\numberwithin{equation}{section}
\allowdisplaybreaks
\theoremstyle{plain}
\newtheorem{theorem}{Theorem}[section]
\newtheorem{proposition}[theorem]{Proposition}
\newtheorem{lemma}[theorem]{Lemma}
\newtheorem{corollary}[theorem]{Corollary}

\theoremstyle{definition}
\newtheorem{definition}[theorem]{Definition}
\newtheorem{example}[theorem]{Example}
\theoremstyle{remark}
\newtheorem{remark}[theorem]{Remark}
\newcommand{\R}{\mathbb R}
\newcommand{\N}{\mathbb N}
\newcommand{\Z}{\mathbb Z}
\newcommand{\C}{\mathbb C}
\newcommand{\Sch}{\mathcal S}
\newcommand{\HH}{\mathcal H}
\newcommand{\supp}{\operatorname{supp}}
\newcommand{\diam}{\operatorname{diam}}

\newcommand{\spanop}{\operatorname{span}}
\newcommand{\dimH}{\dim_{\mathrm H}}

\newcommand{\cov}{\mathcal N}
\newcommand{\la}{\langle}
\newcommand{\ra}{\rangle}
\newcommand{\norm}[1]{\left\lVert#1\right\rVert}
\newcommand{\abs}[1]{\left\lvert#1\right\rvert}

\newcommand{\one}{\mathbf 1}
\newcommand{\eps}{\varepsilon}

\newcommand{\dd}{\,d}
\newcommand{\ee}[1]{e^{2\pi i #1}}
\title[Uniqueness of functions in weak $L^p$]
{Weak $L^p$ integrability of functions with Fourier support in fractal subsets of finite-type manifolds}
\author{Thanh-Tinh Cao}
\address{Department of Mathematics and Physics, University of Information Technology, Ho Chi Minh City, Vietnam,
and Vietnam National University, Ho Chi Minh City, Vietnam}
\email{tinhct@uit.edu.vn}
\author{Minh-Quy Pham}
\address{Department of Mathematics, Virginia Tech,
Blacksburg, VA 24061-1026, USA}
\email{mqpham@vt.edu}
\date{\today}
\subjclass[2020]{Primary 42B10; Secondary 28A78, 46E30}
\keywords{Fourier support, weak $L^p$ spaces, Hausdorff dimension,
finite-type manifolds, Fourier uniqueness}

\begin{document}

\begin{abstract}
Let $d\geq 2$.
We prove that for a finite-type submanifold $\mathcal M\subset\R^d$ with a common rank profile, there exists an exponent $P_{\mathcal M}$, determined by the geometry of $\mathcal M$, with the following property.
If $E\subset\mathcal M$ is compact with Hausdorff dimension $\alpha>0$  and $f\in L^{p,\infty}(\R^d)$ represents a tempered distribution whose Fourier transform $\widehat{f}$ is supported in $E$, then $f=0$ for $1\le p<1+P_{\mathcal M}/\alpha$.
Equivalently, a nonzero $f$ with these support and integrability properties satisfies the uncertainty principle
\[
 (p-1)\dimH E\ge P_{\mathcal M}.
\]
On a fixed manifold, decreasing $p$ raises the lower bound on the Hausdorff dimension of the Fourier support.
Conversely, a smaller support dimension requires a larger integrability exponent $p$.
If $\HH^\alpha(E)<\infty$, then uniqueness also holds at $p=1+P_{\mathcal M}/\alpha$ in strong $L^p$.

In particular, we recover the uniqueness range of Guo, Iosevich, Zhang, and Zorin-Kranich for the full moment curve.
For compact subsets of the moment curve in $\R^d$ with $d\ge3$, we replace the covering exponent in the result of Duan, Pham, and Ryou by Hausdorff dimension.
This yields a larger strict uniqueness range whenever the Hausdorff dimension is strictly smaller than upper box dimension of $E$.

We apply our results to obtain necessary conditions for Fourier integrability of measures and for Fourier restriction estimates.
We also obtain an obstruction for convolution powers: if a probability measure with compact support $E\subset\mathcal M$ has an $\ell$-fold convolution with an $L^2$ density, then $(2\ell-1)\dimH E\ge P_{\mathcal M}$.
\end{abstract}
\maketitle

\section{Introduction and main results}\label{sec:intro}
Let $d\geq1$ be an integer and $1\leq p<\infty$.
How thin can the Fourier support $\supp\widehat f\subset\R^d$ of a nonzero $f\in L^p(\R^d)$ be?
As we will see, the smaller $p$ is, the higher the lower bound on the Hausdorff dimension of the Fourier support.

The above question is equivalent to the following uniqueness question: for which closed sets $E\subset\R^d$ do the Fourier data outside $E$ uniquely determine $f\in L^p(\R^d)$?
Observe that if two functions $f,g\in L^p(\R^d)$ have the same Fourier transform on $\R^d\setminus E$ if and only if $\supp\widehat{(f-g)}\subset E$.
Since $f-g$ also lies in $L^p(\R^d)$, it suffices to ask when the only such difference is zero.
More precisely, given a compact set $E\subset\R^d$ of Hausdorff dimension $\alpha\in(0,d)$, for which exponents $1\le p<\infty$ does
\[
 f\in L^p(\R^d),\qquad \supp\widehat f\subset E
 \quad\Longrightarrow\quad f=0
\]
hold?

The classical Hausdorff dimension bound gives $f=0$ for $1\leq p<2d/\alpha$.
On the line, this follows from Beurling's capacity argument \cite{Beurling51}.
In higher dimensions, Edgar and Rosenblatt \cite{EdgarRosenblatt79} proved uniqueness for $p<2d/(d-1)$ when $\dimH E\leq d-1$.
The general range, in turn, follows from Kahane's capacity identity \cite{Kahane84}.
Additional hypotheses on $E$ can give uniqueness at the endpoint: H\"ormander \cite{Hormander73} treated the hypersurface case, and Agranovsky and Narayanan \cite{AN04} proved uniqueness up to $p=2d/m$ when $E$ lies in an $m$-dimensional $C^1$ submanifold, where $1\leq m<d$.
For fractal sets, Senthil Raani \cite[Theorem~2.3]{Raani14} obtained the endpoint $p=2d/\alpha$ under the assumption of finite $\alpha$-dimensional packing premeasure.
Senthil Raani also used products of Salem sets to construct counterexamples for $p>2d/\alpha$.
At $p=2d/\alpha$, Dobronravov \cite{DobSmall26} constructed a probability measure $\mu$ supported on a compact set $E$ with $\HH^\alpha(E)=0$ and $\check\mu\in L^{2d/\alpha}(\R^d)$.
Thus Hausdorff dimension alone does not give uniqueness at the endpoint.

The general threshold $2d/\alpha$ may not be sharp for Fourier supports on manifolds of high codimension.
For the moment curve $\Gamma_d=\{(t,t^2,\ldots,t^d):0\leq t\leq1\}$, where $d\geq3$, Guo, Iosevich, Zhang, and Zorin-Kranich \cite{GIZZK23} improved the general range $1\leq p\leq2d$ to $1\leq p<(d^2+d+2)/2$.
Duan, Pham, and Ryou \cite{DPR25} extended this result to compact sets $E\subset\Gamma_d$ that can be covered by $O(\delta^{-\alpha})$ balls of radius $\delta$ for every $0<\delta<1$, where $0<\alpha\le1$.
They proved that if $f\in L^p(\R^d)$ and $\supp\widehat f\subset E$, then $f=0$ provided $1\le p\le\frac{d^2+d+2\alpha}{2\alpha}$.
Their random Cantor set construction also shows that this range is sharp when the Hausdorff dimension $\dimH E$ equals the covering exponent $\alpha$.

However, neither the moment-curve results nor the manifold results discussed above give a uniqueness threshold in terms of the Hausdorff dimension of $E$.
This matters because Hausdorff dimension can be strictly smaller than the covering exponent or the upper box dimension.
For example, fix $0<\varepsilon<1$ and take a compact set $E\subset\Gamma_d$ with Hausdorff dimension $\dimH E=\varepsilon$ and upper box dimension $\overline{\dim_M}E=1$.
The theorem of Duan, Pham, and Ryou \cite{DPR25} applies with covering exponent one, but not with a smaller exponent.
Therefore, it forces $f\in L^p(\R^d)$ with $\supp\widehat f\subset E$ to vanish for $1\le p\le 1+d(d+1)/2$, regardless of how small $\varepsilon$ is.
When $\varepsilon$ is small, one might expect uniqueness for a much wider range of $p$.

To address this gap, we study the range of integrability exponents $p$ for which a function $f$ in the space $L^{p,\infty}(\R^d)$, $1\leq p<\infty$, must vanish when its Fourier support is contained in a compact set $E$ of Hausdorff dimension $\alpha$ on a finite-type manifold.
In particular, for the compact set $E\subset\Gamma_d$ above, we show that if $f\in L^{p,\infty}(\R^d)$ represents a tempered distribution and $\supp\widehat f\subset E$, then $f=0$ for
\[
 1\le p<1+\frac{d(d+1)}{2\varepsilon}.
\]
More generally, the two strict thresholds agree when $\dimH E=\overline{\dim_M}E=\alpha>0$.
At the strong $L^p$ endpoint, we replace the uniform covering bound at exponent $\alpha$ by the weaker condition $\HH^\alpha(E)<\infty$.
For the full moment curve, $\alpha=1$, and our strict range agrees with that in \cite{GIZZK23}.

Before stating our general result, we recall some terminology.
Let $d\geq 2$, and let $\mathcal M\subset\R^d$ be an $m$-dimensional submanifold of class $C^k$, where $1\le m<d$ and $k\geq 1$.
\begin{definition}
Fix $x_0\in\mathcal M$.
Let $\varphi:U\to\mathcal M$ be a local $C^k$ parametrization of a neighborhood of $x_0$, where $U\subset\R^m$ is open, and let $u_0\in U$ satisfy $\varphi(u_0)=x_0$.
For each $1\le j\le k$ we write
\[
 V^{(j)}(x_0)=\spanop
 \bigl\{\partial^\nu\varphi(u_0):\nu\in\mathbb N_0^m,\ 1\le |\nu|\le j\bigr\}
\]
for the $j$th-order osculating space of $\mathcal M$ at $x_0$, and set $r_j(x_0)=\dim V^{(j)}(x_0)$.

We call $(r_1(x_0),\ldots,r_k(x_0))$ the \emph{rank profile} of $\mathcal M$ at $x_0$.
If this sequence is the same at every point, we call it the rank profile of $\mathcal M$.
\end{definition}
Note that the spaces $V^{(j)}(x_0)$ are independent of the parametrization, by the chain rule applied to the change of coordinates and its inverse.

We say that $\mathcal M$ is of \emph{finite type at $x_0$} if $r_j(x_0)=d$ for some $1\le j\le k$, and call the least such $j$ the \emph{type} at $x_0$.
The manifold is of \emph{finite type} if it is of finite type at every point.

Throughout, assume that $\mathcal M$ has a common rank profile $(r_1,\ldots,r_k)$ with $r_{k-1}<d=r_k$. Thus $\mathcal M$ is of type $k$ at every point. With $r_0=0$, we define the \emph{volume exponent} of $\mathcal M$ by \[
 P_{\mathcal M}=\sum_{j=1}^k j(r_j-r_{j-1}).
\]
Since $\mathcal M$ is $m$-dimensional, its first derivatives span the tangent space, so $r_1=m$.
The remaining $d-m$ directions enter with weights at least two, so $P_{\mathcal M}\ge 2d-m>d$.

Our first main result is the following.

\begin{theorem}\label{thm:hausdorff}
Let $d,m,k$ be integers with $1\le m<d$ and $k\ge2$.
Let $1\le p<\infty$.
Let $\mathcal M\subset\R^d$ be an $m$-dimensional $C^k$ submanifold of type $k$ at every point, with rank profile $(r_1,\ldots,r_k)$, where $m=r_1\le\cdots\le r_{k-1}<r_k=d$.
Let $E\subset\mathcal M$ be compact with $\dimH E=\alpha\in (0,m]$.

Let $f\in L^{p,\infty}(\R^d)$ represent a tempered distribution and suppose that $\supp\widehat f\subset E$.
If
\begin{equation*}
p<p_{\mathcal{M},\alpha}=1+\frac{P_{\mathcal M}}{\alpha},
\end{equation*}
then $f=0$.
\end{theorem}
Recent work has improved the $L^p$ uniqueness range for functions with Fourier supports on submanifolds of high codimension.
Cao, Li, Miao, and Pang \cite[Proposition~7.13]{CLMP25} give a dimension-based uniqueness bound for smooth graph patches.
Deodhar \cite{Deo26} extends the $L^p$ uniqueness range beyond the classical threshold $2d/m$ when $E$ is contained in an $m$-dimensional $C^d$ submanifold of $\R^d$ with $m<d/2$, and Bi and Sun \cite{BiSun26} obtain improved ranges on Parsell--Vinogradov manifolds.
Both results apply to fractal subsets of these manifolds, but their ranges do not depend on the Hausdorff dimension of the subset.
Dobronravov \cite{Dob26}, on the other hand, proves uniqueness without a manifold hypothesis at $p=2d/\alpha>2$ in $L^{p,q}$ for $1\le q<\infty$ when the compact Fourier support is $\sigma$-finite with respect to the Netrusov--Hausdorff capacity.
In a different direction, Deodhar and Iosevich \cite{DeodharIosevich26} study Fourier transforms of measures with a covering bound and an additional spectral complexity parameter, the Fourier ratio.

Theorem~\ref{thm:hausdorff} combines the Hausdorff dimension $\alpha$ of the Fourier support of the function with the finite-type geometry of the containing manifold, encoded by $P_{\mathcal M}$.
The threshold $p_{\mathcal M,\alpha}$ exceeds the classical Hausdorff dimension threshold $2d/\alpha$ when $P_{\mathcal M}+\alpha>2d$.
It also exceeds a full-manifold threshold $P>1$ from \cite{CLMP25,Deo26,BiSun26} when $\alpha<P_{\mathcal M}/(P-1)$.
The range can improve even when $E$ is an entire surface patch.
The following two-dimensional example in $\R^{14}$ shows how the rank profile yields this improvement.
\begin{example}
Let $I=[1/4,1/2]$ and define $F_D:\R^2\to\R^{14}$ by
\[
 F_D(u,v)=\bigl(u,u^2,u^3,u^4,u^6,\ v,v^2,v^3,v^4,v^6,\ uv,u^3v,uv^3,u^3v^3\bigr).
\]
Set $\mathcal M=F_D((0,1)^2)$ and take $E=F_D(I^2)$ to be the entire compact patch.
After permuting the coordinates, $\mathcal M$ is a smooth graph over $(0,1)^2$, so $\dimH E=2$.
The surface has rank profile $(2,5,9,12,14)$, so $P_{\mathcal M}=42$; see Section~\ref{sec:whole-surface-comparison}.
Theorem~\ref{thm:hausdorff} therefore gives uniqueness in weak $L^p$ for $1\le p<1+42/2=22$.
Since $F_D$ is Lipschitz on $I^2$, $\HH^2(E)<\infty$.
Therefore Theorem~\ref{thm:finite-hausdorff} gives strong-$L^p$ uniqueness also at $p=22$.
On the same patch, Cao, Li, Miao, and Pang \cite{CLMP25} give $p<21$, and Deodhar's general bound \cite[Theorem~1.1]{Deo26} gives $p\le37/2$.
Bi and Sun's theorem on the containing quadratic Parsell--Vinogradov manifold \cite{BiSun26} gives $p<10$.
The gain over the dimension-based bound comes from two directions that first appear at order five.
\end{example}

A classical formulation of the uncertainty principle is as follows.
If a nonzero function $f$ is concentrated on a small scale in physical space, then its Fourier transform $\widehat f$ must spread over a large scale in frequency space, and vice versa.
Theorem~\ref{thm:hausdorff} gives an uncertainty principle for distributions supported on finite-type manifolds.
It relates the Hausdorff dimension of the support of the distribution to the exponent $p$ for which the Fourier transform belongs to weak $L^p$.
On a fixed manifold, a smaller support dimension raises the lower bound on $p$, and vice versa.

\begin{corollary}[An uncertainty principle]\label{cor:uncertainty}
Let $d,m,k$ be integers with $1\le m<d$ and $k\ge2$.
Let $\mathcal M\subset\R^d$ be an $m$-dimensional $C^k$ submanifold of type $k$ at every point, with common rank profile $(r_1,\ldots,r_k)$, where $m=r_1\le\cdots\le r_{k-1}<r_k=d$.
Let $\mu$ be a nonzero distribution on $\R^d$ with compact support $E\subset\mathcal M$ and $0<\dimH E\le m$.
If $\widehat\mu\in L^{p,\infty}(\R^d)$ for some $1\le p<\infty$, then
\[
 (p-1)\dimH E\ge P_{\mathcal M}.
\]
\end{corollary}
\begin{proof}
To apply the theorem, note that the function $f=\check\mu=\widehat\mu(-\,\cdot)$ belongs to $L^{p,\infty}(\R^d)$ and has Fourier transform $\mu$.
Since $f\ne0$, Theorem~\ref{thm:hausdorff} gives the inequality.
This completes the proof of the corollary.
\end{proof}
For strong $L^p$ functions, finite $\alpha$-dimensional Hausdorff measure of $E$ suffices to reach the endpoint of Theorem~\ref{thm:hausdorff}.
No bound on the covering numbers of $E$ is needed.
\begin{samepage}
    \begin{theorem}\label{thm:finite-hausdorff}
Let $d,m,k$ be integers with $1\le m<d$ and $k\ge2$.
Let $\mathcal M\subset\R^d$ be an $m$-dimensional $C^k$ submanifold of type $k$ at every point, with common rank profile $(r_1,\ldots,r_k)$, where $m=r_1\le\cdots\le r_{k-1}<r_k=d$.
Let $E\subset\mathcal M$ be compact and suppose that $\HH^\alpha(E)<\infty$ for some $0<\alpha\le m$.
If $f\in L^p(\R^d)$, $\supp\widehat f\subset E$, and
\[
 1\le p\le p_{\mathcal M,\alpha}=1+\frac{P_{\mathcal M}}{\alpha},
\]
then $f=0$.
\end{theorem}
\end{samepage}

At the weak-$L^p$ endpoint, our argument gives uniqueness if the normalized covering number tends to zero along a sequence of scales.
We write $\cov(E,\delta)$ for the least number of open balls of radius $\delta$ needed to cover $E$.

\begin{samepage}
\begin{theorem}\label{thm:covering}
Let $d,m,k$ be integers with $1\le m<d$ and $k\ge2$.
Let $1\le p<\infty$.
Let $\mathcal M\subset\R^d$ be an $m$-dimensional $C^k$ submanifold of type $k$ at every point, with rank profile $(r_1,\ldots,r_k)$, where $m=r_1\le\cdots\le r_{k-1}<r_k=d$.
Suppose that $E\subset\mathcal M$ is compact. \\ Let $f\in L^{p,\infty}(\R^d)$ represent a tempered distribution and suppose that $\supp\widehat f\subset E$.
Fix $0<\alpha\le m$.
Then $f=0$ in either of the following cases:
\begin{enumerate}
\item[(i)] $1\le p<p_{\mathcal M,\alpha}$ and
\begin{equation}\label{eq:lower-cover-assumption}
 \liminf_{\delta\downarrow0}\delta^\alpha\cov(E,\delta)<\infty.
\end{equation}
\item[(ii)] $p=p_{\mathcal M,\alpha}$ and
\begin{equation}\label{eq:zero-content}
 \liminf_{\delta\downarrow0}\delta^\alpha\cov(E,\delta)=0.
\end{equation}
\end{enumerate}
In particular, if $\alpha=\dimH E>0$ and \eqref{eq:zero-content} holds, then Theorem~\ref{thm:hausdorff} holds for $p=p_{\mathcal M,\alpha}$.
\end{theorem}
\end{samepage}
A uniform covering bound of the form
\begin{equation}\label{eq:cover-assumption}
 \cov(E,\delta)\le C_E\delta^{-\alpha}
 \qquad(0<\delta<\delta_0),
\end{equation}
for some $C_E,\delta_0>0$, implies \eqref{eq:lower-cover-assumption}.
Theorem~\ref{thm:covering} needs less: for $p<p_{\mathcal M,\alpha}$, it is enough that $\delta^\alpha\cov(E,\delta)$ stay bounded along a sequence $\delta\downarrow0$.
If this quantity tends to zero along such a sequence, the theorem also gives uniqueness at $p=p_{\mathcal M,\alpha}$.

We use these uniqueness results in Section~\ref{sec:applications} to give necessary conditions for Fourier integrability of measures, and necessary conditions for Fourier extension estimates to hold, and determine when a prescribed convolution power cannot have an $L^2$ density.

The range in Theorem~\ref{thm:hausdorff} cannot be enlarged in general.
On the moment curve in $\R^d$ with $d\ge3$, Duan, Pham, and Ryou \cite[Theorems~1.1--1.2 and Proposition~3.1]{DPR25} establish sharpness from above for supports whose Hausdorff dimension $\alpha$ is also their covering exponent.
More precisely, their examples satisfy $\cov(E,\delta)\approx\delta^{-\alpha}$, and for every $p>1+d(d+1)/(2\alpha)$ they give a nonzero $L^p$ function with Fourier support in such an $E$.
Since $L^p\subset L^{p,\infty}$, these examples also rule out any larger weak-$L^p$ range.

In Section~\ref{sec:product-sharpness}, we give a two-dimensional example by taking the product of two cubic moment curves.
The product has $P_{\mathcal M}=12$, and for every $0<\alpha\le2$ and $p>1+12/\alpha$ we construct a compact set $E$ on it with $\dimH E=\alpha$ that supports the Fourier transform of a nonzero weak-$L^p$ function.

\subsection*{Proof outline}\label{sec:proof-outline}

We cover the part of $\mathcal M$ containing $E$ by anisotropic boxes $Q$ adapted to the geometry of $\mathcal M$.
The main task is to choose $\delta_0>0$ and construct these boxes and their cutoffs $b_Q$ at scales $0<\delta\le\delta_0$ with the following properties:
\begin{enumerate}
\item[(i)] Each box $Q$ must be centered on $\mathcal M$ and contain the corresponding portion of the graph at scale $\delta$.
Its side lengths are comparable to $\delta^j$ in the $r_j-r_{j-1}$ directions that first enter the osculating spaces at order $j$, for every $1\le j\le k$, with $r_0=0$.
These lengths accommodate the Taylor displacement in each direction and give the volume factor $|Q|\approx\delta^{P_{\mathcal M}}$ used in the kernel estimates.
Lemma~\ref{lem:taylor} supplies these boxes.
\item[(ii)] Let $b_Q$ be the smooth cutoff associated with $Q$ by \eqref{eq:bI}.
In both proofs we choose nonnegative parameter cutoffs, so $b_Q\ge0$ and $b_Q=0$ outside $CQ$, for $C>0$.
We need the bounds
\[
 \|\check b_Q\|_{L^1}\lesssim1,\qquad
 \|b_Q\|_{L^\infty}\lesssim1,\qquad
 |\supp b_Q|\lesssim\delta^{P_{\mathcal M}},
\]
with constants independent of the center and the scale.
The first controls the $L^1$ contribution of each box; the other two give $\|\check b_Q\|_{L^\infty}\lesssim\delta^{P_{\mathcal M}}$ by Fourier inversion.
These are the single-cutoff bounds used to estimate the sum.
They are proved in Lemma~\ref{lem:symbol}.
\item[(iii)] We need
\[
 \chi:=\sum_Q b_Q=1\quad\text{on an open neighborhood of }E.
\]
Equality near $E$, not just on $E$, gives $\chi\widehat f=\widehat f$ because $\widehat f$ is a distribution supported in $E$.
Thus $f=f*K$, where $K=\check\chi$.
Lemma~\ref{lem:cutoff-neighborhood} gives the neighborhood equality when the cutoffs in parameter space sum to one near the preimage of $E$.
\end{enumerate}

\smallskip
\noindent\textit{Theorem~\ref{thm:hausdorff}.} First we pass from weak $L^p$ to strong $L^{p_0}$, choosing $p_0>p$ so that $\dimH E<s:=P_{\mathcal M}/(p_0-1)<P_{\mathcal M}$.
Since $\HH^s(E)=0$, we can cover $E$ by boxes $Q$ with scales $0<\delta_Q\le\delta_*\le\delta_0$ and $\sum_Q\delta_Q^s\le1$.
This is the role of Lemma~\ref{lem:dyadic-partition} in the construction.

The boxes need not have the same thickness.
To ensure that $\chi=\sum_Q b_Q$ equals one on an open neighborhood of $E$, we use the smallest box scale to choose a common transverse thickness $c\min_Q\delta_Q^k>0$.
This is provided by Lemma~\ref{lem:cutoff-neighborhood}.
Thus $f=f*K$, where $K=\check\chi$.
Since $P_{\mathcal M}(p_0'-1)=s$, the kernel satisfies
\[
 \|K\|_{L^{p_0'}}^{p_0'}\lesssim\sum_Q\delta_Q^s,
 \qquad
 \|K\|_{L^\infty}\lesssim\delta_*^{P_{\mathcal M}-s}\sum_Q\delta_Q^s,
\]
by Lemma~\ref{lem:all-scale-assembly}.
The first bound is uniform and the second tends to zero as $\delta_*\downarrow0$.
Hence $f=0$ by Lemma~\ref{lem:cutoff}(ii).

\smallskip
\noindent\textit{Theorem~\ref{thm:covering}.} The fixed-scale argument follows the cutoff and interpolation method of \cite[proof of Proposition~7.13]{CLMP25}.
For $2<p<p_{\mathcal M,\alpha}$, the assumption $\liminf_{\delta\downarrow0}\delta^\alpha\cov(E,\delta)<\infty$ gives a sequence $\delta_n\downarrow0$ along which the normalized covering counts stay bounded.
For $p>2$, we take all boxes at scale $\delta$.
Their associated cutoffs have bounded overlap, and the number of boxes is at most a constant times $\cov(E,\delta)$.
Write $\chi_\delta=\sum_Q b_Q$ and $K_\delta=\check\chi_\delta$.
Again $\chi_\delta=1$ near $E$, but now we need $\|K_\delta\|_{L^{p',1}}\to0$ along a sequence.
Setting $A_\delta=\delta^\alpha\cov(E,\delta)$ gives
\[
 \|K_\delta\|_{L^{p',1}}\lesssim
 A_\delta^{1-1/p}\delta^{[P_{\mathcal M}-\alpha(p-1)]/p}
\]
by Lemma~\ref{lem:assembly}.
Below $p_{\mathcal M,\alpha}$ the power of $\delta$ is positive, so bounded $A_\delta$ along a sequence $\delta\downarrow0$ suffices.
At the endpoint this power is zero, so we instead take a sequence with $A_\delta\to0$.
In both cases $f=0$ by Lemma~\ref{lem:cutoff}(i).
For $1\le p\le2$, the Hausdorff argument applies without a covering-number assumption.

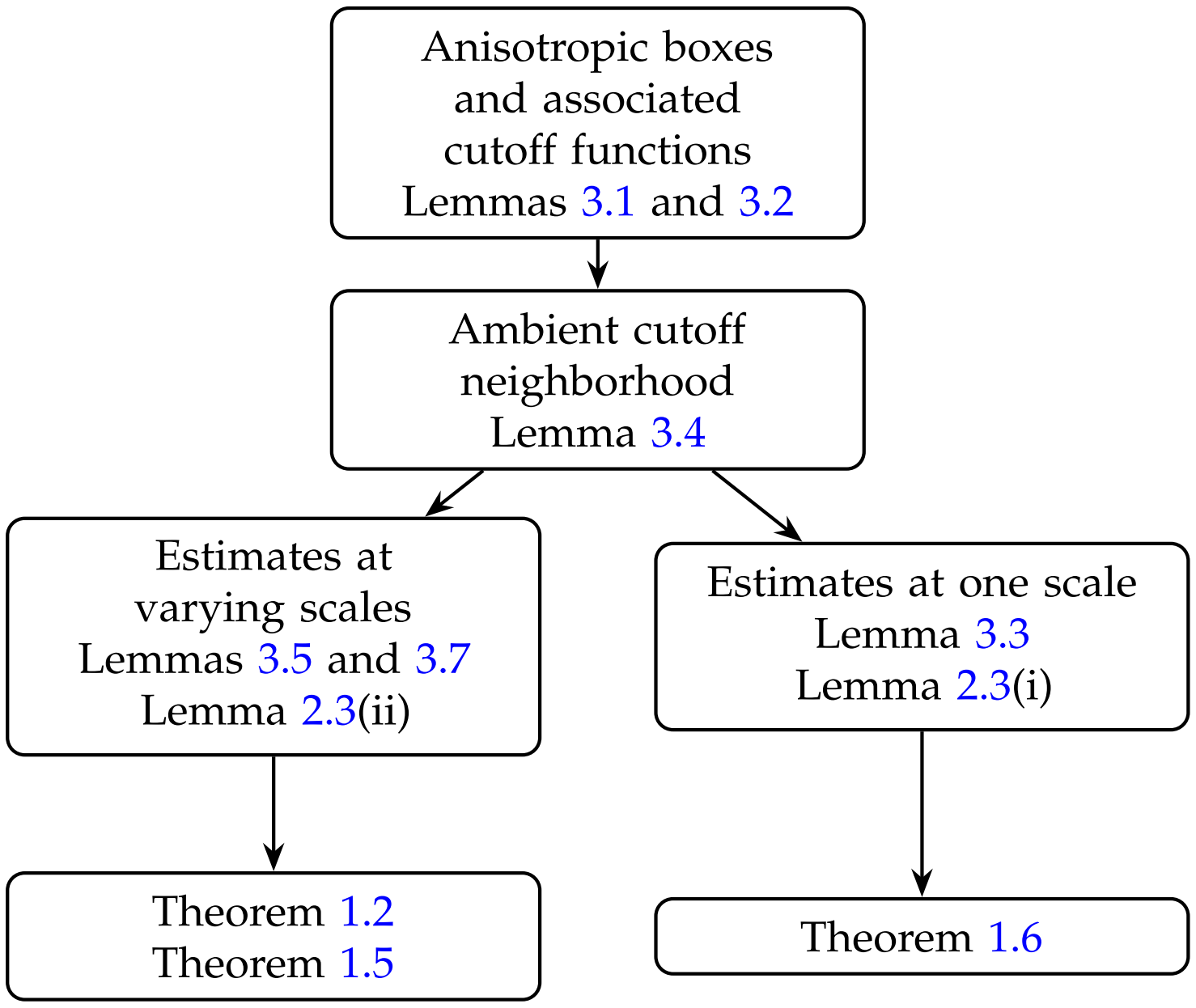
\begin{figure}[htbp]
\centering
\begin{tikzpicture}[>=Stealth,box/.style={draw,rounded corners,align=center,
 font=\small,text width=4.0cm,inner sep=5pt},every path/.style={thick}]
\node[box] (geometry) at (0,0) {Anisotropic boxes and associated cutoff functions\\
Lemmas \ref{lem:taylor} and \ref{lem:symbol}};
\node[box] (assembly) at (0,-2.1) {Ambient cutoff neighborhood\\Lemma \ref{lem:cutoff-neighborhood}};
\node[box] (haus) at (-2.65,-4.2) {Estimates at varying scales \\
Lemmas \ref{lem:dyadic-partition} and \ref{lem:all-scale-assembly}\\
Lemma \ref{lem:cutoff}(ii)};
\node[box] (cover) at (2.65,-4.2) {Estimates at one scale\\
Lemma \ref{lem:assembly}\\
Lemma \ref{lem:cutoff}(i)};
\node[box] (ht) at (-2.65,-6.65) {Theorem \ref{thm:hausdorff}\\Theorem \ref{thm:finite-hausdorff}};
\node[box] (ct) at (2.65,-6.65) {Theorem \ref{thm:covering}};
\draw[->] (geometry)--(assembly);
\draw[->] (assembly)--(haus);
\draw[->] (assembly)--(cover);
\draw[->] (haus)--(ht);
\draw[->] (cover)--(ct);
\end{tikzpicture}
\caption{The two uniqueness arguments share the kernel construction,
but require different information about the spectrum. The Lorentz
estimates of Lemmas \ref{lem:lorentz} and \ref{lem:interpolation} enter the construction
and the cutoff criterion. Arrows indicate the flow of the two proofs,
not implications between the lemmas grouped in a box.}
\label{fig:proof-architecture}
\end{figure}

\subsection*{Organization} In Section~\ref{sec:preliminaries}, we introduce the notation and record the geometric and Lorentz-space estimates used in the proofs.
In Section~\ref{sec:proof-main}, we construct cutoff functions at varying scales and prove Theorems~\ref{thm:hausdorff} and~\ref{thm:finite-hausdorff}, including the strong $L^p$ endpoint.
Then Section~\ref{sec:proof-covering} uses a single scale at a time to prove Theorem~\ref{thm:covering}.
Finally, Section~\ref{sec:applications} gives applications, and Section~\ref{sec:examples} presents a sharpness example and a comparison with earlier uniqueness bounds.

\newpage
\section{Preliminaries}\label{sec:preliminaries}

 We first recall the terminology for Hausdorff dimension, covering numbers, and finite-type manifolds.
 We then record the facts about weak $L^p$ that we use to estimate the cutoff kernels.

\subsection{Notation}\label{sec:notation}

Let $n\ge1$ be an integer and let $U\subset\R^n$ be open.
We write $C_c^\infty(U)$ for the space of smooth functions $f:U\to\C$ whose support is a compact subset of $U$.

We write $x,\xi\in\R^d$ for the physical and frequency variables, respectively, and $|A|=\mathcal L^d(A)$ for the Lebesgue measure of a measurable set $A\subset\R^d$.
If $A$ and $B$ are nonnegative, we write $A\lesssim B$ when $A\le CB$ and $A\approx B$ if both $A\lesssim B$ and $B\lesssim A$.
The constant $C$ may depend on the fixed geometric data and stated exponents, but not on the covering scale or the number of cutoff pieces.
The H\"older conjugate $p'$ of $p\in[1,\infty]$ satisfies $1/p+1/p'=1$, with $1/\infty=0$.

For $1\le p<\infty$, the Lebesgue space $L^p(\R^d)$ consists of the measurable functions $f:\R^d\to\C$ such that
\[
 \int_{\R^d}|f(x)|^p\,\dd x<\infty.
\]

The weak Lebesgue space $L^{p,\infty}(\R^d)$, for $1\le p<\infty$, consists of measurable functions $f$ on $\R^d$, modulo equality almost everywhere, with finite quasi-norm
\begin{equation}\label{eq:weak-lebesgue}
 \|f\|_{p,\infty}:=\sup_{\lambda>0}
 \lambda\,\bigl|\{x\in\R^d:|f(x)|>\lambda\}\bigr|^{1/p}.
\end{equation}
Our conventions for the Fourier and inverse Fourier transforms of $f,g\in L^1(\R^d)$ are
\[
 \widehat f(\xi)=\int_{\R^d}e^{-2\pi i x\cdot\xi}f(x)\dd x,
 \qquad
 \check g(x)=\int_{\R^d}\ee{x\cdot\xi}g(\xi)\dd\xi.
\]
Write $\Sch(\R^d)$ for the Schwartz space and $\Sch'(\R^d)$ for the space of tempered distributions.
Given $T\in\Sch'(\R^d)$, the Fourier transform of $T$ is defined by the pairing $\la\widehat T,\varphi\ra=\la T,\widehat\varphi\ra$ for every $\varphi\in\Sch(\R^d)$.
The inverse Fourier transform is defined in the same way with $\check\varphi$ in place of $\widehat\varphi$.
With this normalization, Plancherel's theorem gives $\|\widehat f\|_2=\|f\|_2$ for $f\in L^2(\R^d)$.
If $E\subset\R^d$ is closed and $f\in\Sch'(\R^d)$, the inclusion $\supp\widehat f\subset E$ means that $\la\widehat f,\varphi\ra=0$ for every $\varphi\in C_c^\infty(\R^d\setminus E)$.
At $p=1$ we require that $f$ be locally integrable and that its induced distribution be tempered; weak $L^1$ alone does not ensure local integrability (see, for example, \cite{Grafakos14}).

\subsection{Dimensions, measures, and covering numbers}\label{sec:dimension-prelim} Let $d\geq 1$, and let $E\subset\R^d$.
For $ s\ge0$, the $s$-dimensional Hausdorff measure of $E$ is $\HH^s(E)=\lim_{\delta\downarrow0}\HH_{\delta}^s(E)$, where
\[
 \HH_{\delta}^s(E)=\inf\left\{\sum_i(\diam U_i)^s:
 E\subset\bigcup_i U_i,\ \diam U_i\le \delta\text{ for every }i\right\}, \quad \delta>0,
\]
where the infimum is over finite or countable families of nonempty sets $U_i\subset\R^d$ and $\diam U_i=\sup_{x,y\in U_i}|x-y|$.
At $s=0$, each nonempty covering set contributes one.
For $s>d$, we have $\HH^s(E)=0$.
In terms of these measures, the Hausdorff dimension $\dimH E\in[0,d]$ is given by
\[
 \dimH E=\inf\{s>0:\HH^s(E)=0\}.
\]
See \cite[Section~2.2]{Mattila15}.

A compact set $E\subset\R^d$ with $\HH^s(E)<\infty$ admits, for every $\delta>0$, a finite cover by open balls $B(x_i,r_i)$ with $x_i\in E$, $0<r_i\le \delta$, and $\sum_i r_i^s\lesssim_s 1+\HH^s(E)$.
If $F:E\to\R^n$ is Lipschitz with constant $C_F$, then $\HH^s(F(E))\le C_F^s\HH^s(E)$, so bi-Lipschitz maps preserve Hausdorff dimension \cite[Sections~3.1--3.2]{Falconer14}.

Let $E\subset\R^d$ be nonempty and bounded.
For $\delta>0$, the \emph{covering number} $\cov(E,\delta)$ is the least number of open balls of radius $\delta$, centered anywhere in $\R^d$, that cover $E$.
The upper Minkowski dimension (also called the upper box dimension) is defined by
\[
 \overline{\dim_M}E=\limsup_{\delta\downarrow0}
 \frac{\log\cov(E,\delta)}{\log(1/\delta)}.
\]
The upper Minkowski dimension lies in $[0,d]$ and is at most $m$ for compact subsets of an $m$-dimensional submanifold; see \cite[Section~2.3]{Mattila15}.

Let $0<\alpha\le d$.
We say that $\alpha$ is a \emph{covering exponent} for $E$ if there is a constant $C>0$ such that $\cov(E,\delta)\le C\delta^{-\alpha}$ for every $0<\delta<1$.
We call $\liminf_{\delta\downarrow0}\delta^\alpha\cov(E,\delta)$ the \emph{lower $\alpha$-covering content}.
It takes values in $[0,\infty]$.

\subsection{Finite-type submanifolds}\label{sec:finite-type}

Let $d,m,k$ be integers with $d\ge2$, $1\le m<d$, and $k\ge2$, and let $\mathcal M\subset\R^d$ be an $m$-dimensional $C^k$ submanifold.
Without loss of generality, we may assume that $\mathcal M=F(U)$ for a single $C^k$ parametrization $F:U\to\mathcal M$, where $U\subset\R^m$ is open.
Indeed, a localization argument reduces the general case to a single coordinate patch.

Fix $x_0\in\mathcal M$ and take orthonormal graph coordinates centered at $x_0$, with $0\in U$ and $F(0)=x_0$:
\begin{equation*}
 F(u)=(u,\Phi(u)),\qquad u\in U,\quad \Phi\in C^k(U;\R^{d-m}),
 \quad \Phi(0)=0.
\end{equation*}
This choice is possible by the inverse function theorem.
The inverse of $F$ on $F(U)$ is the restriction of the coordinate projection $\pi:\R^d\to\R^m$ onto the first $m$ coordinates.
If $B\subset U$ is a closed ball, the mean value estimate gives
\[
 |u-v|\le |F(u)-F(v)|\le C_B|u-v|\quad\text{for all }u,v\in B,
 \qquad C_B=\bigl(1+\sup_B\|D\Phi\|^2\bigr)^{1/2}.
\]
Every compact set $E\subset\mathcal M$ satisfies $\cov(E,\delta)\le C_E\delta^{-m}$ for $0<\delta<1$ (cf.~\cite[Section~5.3]{Mattila95}).
Thus \eqref{eq:cover-assumption} holds with $\alpha=m$.

Let $u\in U$.
For $1\le j\le k$, the \emph{$j$th-order osculating space} of $\mathcal M$ at $F(u)$ is the linear subspace
\begin{equation*}
 V^{(j)}(F(u))=\spanop\{\partial^\nu F(u):\nu\in\N_0^m,\ 1\le|\nu|\le j\},
\end{equation*}
where $\N_0=\{0,1,2,\ldots\}$.
This is the notation of \cite[(4.7)]{GuoZhang19} for Parsell--Vinogradov manifolds.

We assume that $\mathcal M$ has the same rank profile $(r_1,\ldots,r_k)$ at every point.
In graph coordinates, this means that
\begin{equation*}
 \begin{gathered}
  \dim V^{(j)}(F(u))=r_j\quad\text{for every }1\le j\le k\text{ and }u\in U,\\
  r_0=0,\qquad m=r_1\le\cdots\le r_{k-1}<r_k=d.
 \end{gathered}
\end{equation*}
Here $k$ is the first order at which the derivatives span $\R^d$.
The common rank profile does more than fix the type: it fixes the intermediate ranks, which can otherwise vary with the point.
Invertible affine changes of ambient coordinates preserve the profile because their linear parts carry each derivative span to the corresponding span.

For each $1\le j\le k$, let $h_j=r_j-r_{j-1}$ be the rank increment at order $j$.
The \textit{volume exponent} is the weighted sum of these increments:
\begin{equation}\label{eq:Q}
\qquad P_{\mathcal M}=\sum_{j=1}^k jh_j.
\end{equation}
The weights in \eqref{eq:Q} come from the side lengths of the anisotropic boxes in Lemma~\ref{lem:taylor}.
For sufficiently small $\delta>0$, the image of a parameter ball of radius $\delta$ centered in a fixed compact subset of $U$ lies in a box with $h_j$ sides of length comparable to $\delta^j$ for each $1\le j\le k$.
It follows that the box has volume comparable to $\delta^{P_{\mathcal M}}$.
Exactly $m$ directions appear at order one, and the remaining $d-m$ appear at order two or higher.
Consequently,
\begin{equation*}
 P_{\mathcal M}\ge m+2(d-m)=2d-m>m.
\end{equation*}

\subsection{Weak $L^p$ estimates and smooth cutoffs}\label{sec:lorentz}

To explain how the cutoffs will be used, let $E\subset\R^d$ be compact, and let $\chi\in C_c^\infty(\R^d)$ be equal to one on a neighborhood of $E$.
If $f\in\Sch'(\R^d)$ and $\supp\widehat f\subset E$, then $\chi\widehat f=\widehat f$.
Taking inverse Fourier transforms gives $f*K=f$ in $\Sch'(\R^d)$, where $K=\check\chi$.
This identity lets us estimate $f$ through bounds on $K$.
The cutoff must equal one on an open neighborhood of the Fourier support, since $\widehat f$ may be a distribution.

The weak $L^p$ estimates below use the decreasing rearrangement.
More precisely, for a measurable function $f$ on $\R^d$, define
\[
 f^*(t)=\inf\bigl\{\lambda\ge0:
       |\{x\in\R^d:|f(x)|>\lambda\}|\le t\bigr\},\qquad t>0,
\]
with $\inf\varnothing=\infty$.
When finite, $f^*(t)$ is the least threshold that $|f|$ exceeds on a set of measure at most $t$.
The function $f^*$ is nonincreasing and has the same superlevel-set measures as $|f|$.

In terms of the rearrangement, the weak $L^p$ quasi-norm in \eqref{eq:weak-lebesgue} is
\[
 \|f\|_{p,\infty}=\sup_{t>0}t^{1/p}f^*(t).
\]
For $1<r<\infty$ and $1\le q<\infty$, the Lorentz space $L^{r,q}(\R^d)$ consists of the measurable functions with finite rearrangement quasi-norm
\[
 \|f\|_{r,q}=\left(\int_0^\infty
     [t^{1/r}f^*(t)]^q\frac{\dd t}{t}\right)^{1/q}.
\]
The case $q=\infty$ is weak $L^r$, with the supremum above; see \cite{Grafakos14} for background.

\begin{lemma}\label{lem:lorentz}
Let $d\ge 1$ be an integer, let $1<p<\infty$, and set $p'=p/(p-1)$.
All function spaces below are on $\R^d$.
\begin{enumerate}
\item[(i)] \textit{(Lorentz H\"older inequality).} For $f\in L^{p,\infty}$ and $g\in L^{p',1}$,
\begin{equation}\label{eq:lorentz-holder}
 \int|fg|\le C_p\norm{f}_{p,\infty}\norm{g}_{p',1}.
\end{equation}
\item[(ii)] \textit{(Young's convolution inequality for weak $L^p$).} If $f\in L^{p,\infty}$ and $K\in L^1$, then
$$\|f*K\|_{p,\infty}\le C_p\|f\|_{p,\infty}\|K\|_1.$$
\item[(iii)] \textit{(Tempered distributions and compact Fourier support).} Every $f\in L^{p,\infty}$ defines a tempered distribution.
If $\supp\widehat f$ is compact, then $f$ has a bounded continuous representative and belongs to $L^r$ for every $p<r<\infty$.
\end{enumerate}
\end{lemma}
\begin{proof}
Part~(i) is the Lorentz H\"older inequality; see the proof of \cite[Theorem~1.4.16(v)]{Grafakos14}.
Indeed, the Hardy--Littlewood rearrangement inequality bounds $\int|fg|$ by $\int_0^\infty f^*(t)g^*(t)\dd t$.
Now use $f^*(t)\le\norm f_{p,\infty}t^{-1/p}$ and the definition of $\norm g_{p',1}$.

Part~(ii) is Young's inequality for weak $L^p$ \cite[Theorem~1.2.13]{Grafakos14}, taking the strong exponent to be $1$ and the input and output weak exponents to be $p$.

For part~(iii), apply \eqref{eq:lorentz-holder} with $g=\one_A$ to get $\int_A|f|\le C_p\norm{f}_{p,\infty}|A|^{1/p'}$ for measurable $A$ of finite measure.
Summing this bound over dyadic annuli against Schwartz decay shows that $f\in\Sch'$.
If $\supp\widehat f$ is compact, choose $\chi\in C_c^\infty$ equal to one near $\supp\widehat f$ and set $K=\check\chi$.
Then $f=f*K$ in $\Sch'$ by \cite[Proposition~2.3.22(11), p.~131]{Grafakos14}.
The convolution is bounded by \eqref{eq:lorentz-holder} and continuous by \cite[Theorem~2.3.20]{Grafakos14}.
Finally, the layer-cake argument in the proof of \cite[Proposition~1.1.4]{Grafakos14} gives
\[
 \norm{f}_r^r\le \frac{r}{r-p}\norm{f}_{p,\infty}^p
                     \norm{f}_\infty^{r-p},\qquad p<r<\infty.
\]
This proves the lemma.
\end{proof}
At $p=1$, the preceding lemma does not apply.
For compact Fourier support, we use the following case of the Bernstein--Nikolski\u\i{} inequality \cite[Theorem~1(i), p.~92]{BangCong05}.
If $f\in L^{1,\infty}(\R^d)$ represents a tempered distribution and $\supp\widehat f$ is contained in a ball of radius $R>0$ centered at the origin, then, for $1<r<\infty$,
\begin{equation}\label{eq:bandlimited-weak-l1}
 \|f\|_r\le C_{d,r} R^{d(1-1/r)}\|f\|_{1,\infty}.
\end{equation}
Apply the unit-ball estimate to $x\mapsto f(x/R)$ and rescale.

To apply Lorentz H\"older to a cutoff kernel, we need its $L^{p',1}$ norm.
For this purpose, the next lemma bounds it in terms of the $L^1$ and $L^2$ norms.

\begin{lemma}[Lorentz interpolation]\label{lem:interpolation}
Let $d\ge1$ be an integer and $2<p<\infty$, and set $p'=p/(p-1)$.
If $K\in L^1(\R^d)\cap L^2(\R^d)$, then
\begin{equation*}
 \norm{K}_{p',1}
 \le C_p\norm{K}_1^{1-2/p}\norm{K}_2^{2/p},
\end{equation*}
where $C_p$ depends only on $p$.
\end{lemma}
\begin{proof}
    This is the real interpolation estimate for $(L^1,L^2)_{2/p,1}$ \cite[Theorem~5.3.1, p.~113]{BerghLofstrom76}, since $1/p'=(1-2/p)+\tfrac12(2/p)$.
This completes the proof of the lemma.
\end{proof}

The next lemma gives two ways to prove $f=0$ from the cutoff identity.
The first uses decay in $L^{p',1}$; the second uses a uniform $L^{p'}$ bound together with decay in $L^\infty$.

\begin{lemma}[Vanishing by Fourier cutoffs]\label{lem:cutoff}
Let $E\subset\R^d$ be compact, let $1<p<\infty$, and set $p'=p/(p-1)$.
For each integer $n\ge1$, let $\chi_n\in C_c^\infty(\R^d)$ equal one on a neighborhood of $E$ and set $K_n=\check\chi_n$.
Suppose that $f\in\Sch'(\R^d)$ and $\supp\widehat f\subset E$.
\begin{enumerate}
\item[(i)] If $f\in L^{p,\infty}(\R^d)$ and $\norm{K_n}_{p',1}\to0$, then $f=0$.
\item[(ii)] If $f\in L^p(\R^d)$, $\sup_n\norm{K_n}_{p'}<\infty$, and $\norm{K_n}_\infty\to0$, then $f=0$.
\end{enumerate}
\end{lemma}
\begin{proof}
For part~(i), it suffices to estimate $f$ directly using the cutoff identity.
Since $\chi_n=1$ on a neighborhood of $E$ and $\supp\widehat f\subset E$, we have $\chi_n\widehat f=\widehat f$, hence $f=f*K_n$ in $\Sch'(\R^d)$.
For each $x\in\R^d$, $K_n(x-\cdot)$ has the same $L^{p',1}$ norm as $K_n$.
Lorentz H\"older inequality \eqref{eq:lorentz-holder} therefore gives absolute convergence of the convolution and
\[
 |f(x)|=|(f*K_n)(x)|
 \le C_p\norm{f}_{p,\infty}\norm{K_n}_{p',1},
\]
for almost every $x$.
The bound is independent of $x$ and tends to zero as $n\to\infty$, so $f=0$ almost everywhere.

For part~(ii), the same cutoff identity holds, but the $L^{p'}$ norms of $K_n$ are only bounded.
To use the decay of $\norm{K_n}_\infty$, we first approximate $f$ by a function in $L^1$.
Given $\eps>0$, choose $g\in C_c^\infty(\R^d)$ with $\norm{f-g}_p<\eps$, and set $M=\sup_n\norm{K_n}_{p'}$.
H\"older's inequality gives absolute convergence of $f*K_n$.
Writing it as $(f-g)*K_n+g*K_n$, we use the $L^{p'}$ bound for the error and the $L^\infty$ bound for the fixed approximation:
\[
 \norm{f}_\infty
 \le \norm{f-g}_p\norm{K_n}_{p'}+\norm{g}_1\norm{K_n}_\infty
 \le \eps M+\norm{g}_1\norm{K_n}_\infty.
\]
With $g$ fixed, let $n\to\infty$ to obtain $\norm{f}_\infty\le\eps M$.
Letting $\eps\downarrow0$ gives $f=0$ almost everywhere.
This completes the proof of the lemma.
\end{proof}

\subsection{Reduction to one patch}\label{sec:chart-localization}

Let $E\subset\mathcal M$ be compact, and suppose $\supp\widehat f\subset E$.
Cover $E$ by finitely many ambient graph neighborhoods and choose cutoffs $\eta_j\in C_c^\infty(\R^d)$ subordinate to this cover, with $\sum_j\eta_j=1$ on a neighborhood of $E$.
Then
\[
 f=\sum_j f_j,\qquad f_j=f*\check\eta_j,\qquad
 \supp\widehat f_j\subset E_j:=E\cap\supp\eta_j.
\]
Each $E_j$ lies in one graph patch.
Since the patches have the same rank profile, they all have volume exponent $P_{\mathcal M}$.
Thus it suffices to prove that each $f_j$ vanishes.

\section{Proofs of Theorems~\ref{thm:hausdorff} and~\ref{thm:finite-hausdorff}}\label{sec:proof-main}

In this section we first prove Theorem~\ref{thm:hausdorff} and~\ref{thm:finite-hausdorff}.
The difficulty is that Hausdorff dimension gives covers by balls of different radii, but does not give a bound on the number of balls at a single scale.
To handle these varying scales, we construct a smooth cutoff for each piece of the cover using an anisotropic box adapted to the osculating spaces.
An ordered partition of unity lets us combine these cutoffs across scales and bound the inverse Fourier transform in terms of the Hausdorff cost of the cover.
Finally, the resulting kernel estimates give vanishing by Lemma~\ref{lem:cutoff}.

\subsection{Anisotropic boxes and uniform cutoff estimates}\label{sec:anisotropic-boxes} Let $d,m,k$ be integers with $d\ge2$, $1\le m<d$, and $k\ge2$.
Assume that the $m$-dimensional $C^k$ graph $\mathcal M=F(U)\subset\R^d$ has type $k$ at every point and common rank profile $(r_1,\ldots,r_k)$, with $m=r_1\le\cdots\le r_{k-1}<r_k=d$.
Recall that the $j$th-order osculating space $V^{(j)}(F(u))$ is spanned by the partial derivatives of $F$ at $u$ of positive order at most $j$.
Thus $\dim V^{(j)}(F(u))=r_j$ for every $u\in U$ and $1\le j\le k$, and the volume exponent is $P_{\mathcal M}=\sum_{j=1}^k j(r_j-r_{j-1})$, with $r_0=0$.

To build a smooth cutoff equal to one near the Fourier support, we cover the compact set $E\subset F(U)$ by anisotropic boxes.
Each box encloses the image of a parameter ball.
Choose a compact set $\Omega\subset U$ with $F^{-1}(E)\subset\operatorname{int}\Omega$ and fix $C_0>0$.
Take $0<\delta<1$ sufficiently small that $B(u_0,C_0\delta)\subset U$ for every $u_0\in\Omega$; this is possible since $\Omega$ is compact in $U$.
For each $u_0\in\Omega$, we construct a box $Q$ centered at $c_Q=F(u_0)$ to contain $F(B(u_0,C_0\delta))$, and the osculating spaces at $c_Q$ determine its directions and scales: a direction that first appears at order $j$ has side length proportional to $\delta^j$.
More precisely, set $V^{(0)}(c_Q)=\{0\}$ and, for each $1\le j\le k$, choose $e_{r_{j-1}+1},\ldots,e_{r_j}$ to be an orthonormal basis of $V^{(j)}(c_Q)\cap (V^{(j-1)}(c_Q))^\perp$.
If $r_j=r_{j-1}$, that order contributes no basis vectors.
Set $w_\ell=j$ whenever $r_{j-1}<\ell\le r_j$.
Then $w_1=\cdots=w_m=1$ and $\sum_{\ell=1}^d w_\ell=P_{\mathcal M}$.
Write $A_{u_0}$ for the matrix with columns $e_1,\ldots,e_d$, and set
\[
 D_\delta=\operatorname{diag}(\delta^{w_1},\ldots,\delta^{w_d}).
\]
Given $C>0$, define the box by
\[
 Q=Q(u_0,\delta;C)=c_Q+A_{u_0}D_\delta[-C,C]^d.
\]
The side parallel to $e_\ell$ has length $2C\delta^{w_\ell}$, and hence $|Q|=(2C)^d\delta^{P_{\mathcal M}}$.

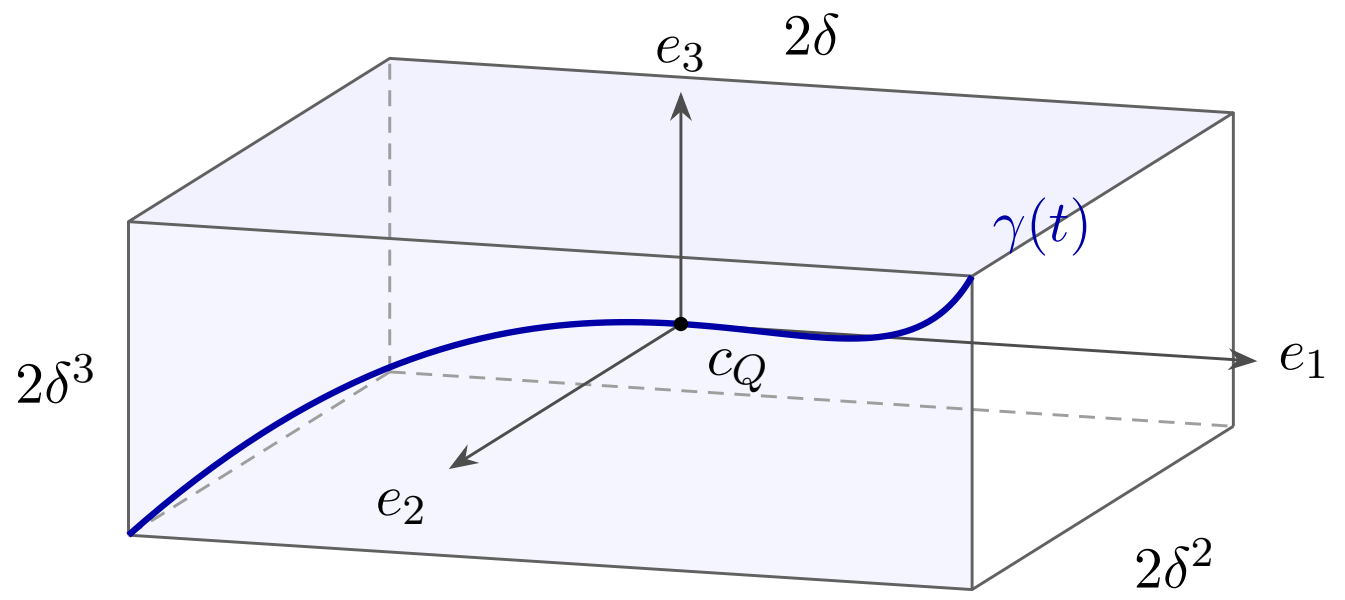
\begin{figure}[htbp]
\centering
\begin{tikzpicture}[
  x={(4.65cm,-0.30cm)}, y={(-2.40cm,-1.50cm)}, z={(0cm,4.80cm)},
  >=Stealth, line join=round,
  boxedge/.style={draw=black!62,line width=.55pt},
  axis/.style={->,draw=black!70,line width=.55pt},
  curve/.style={draw=blue!65!black,line width=1.2pt}]
\coordinate (A) at (-.60,-.36,-.216);
\coordinate (B) at ( .60,-.36,-.216);
\coordinate (C) at ( .60, .36,-.216);
\coordinate (D) at (-.60, .36,-.216);
\coordinate (E) at (-.60,-.36, .216);
\coordinate (F) at ( .60,-.36, .216);
\coordinate (G) at ( .60, .36, .216);
\coordinate (H) at (-.60, .36, .216);
\fill[blue!5] (E)--(F)--(G)--(H)--cycle;
\fill[blue!4] (H)--(G)--(C)--(D)--cycle;
\draw[boxedge,densely dashed,black!38] (A)--(B) (A)--(D) (A)--(E);
\draw[boxedge] (E)--(F)--(G)--(H)--cycle (B)--(C)--(D)
  (B)--(F) (C)--(G) (D)--(H); \path (E)--(F) node[midway,above=2pt] {$2\delta$} (B)--(C) node[midway,below right=2pt] {$2\delta^2$} (D)--(H) node[midway,left=2pt] {$2\delta^3$};
\draw[axis] (0,0,0)--(.82,0,0) node[right] {$e_1$};
\draw[axis] (0,0,0)--(0,.64,0) node[below left] {$e_2$};
\draw[axis] (0,0,0)--(0,0,.32) node[above] {$e_3$};
\draw[curve] plot[domain=-.60:.60,samples=80,variable=\curveparam]
  (\curveparam,{\curveparam*\curveparam},{\curveparam*\curveparam*\curveparam});
\fill[black] (0,0,0) circle[radius=1.35pt] node[below right=1pt] {$c_Q$};
\node[blue!65!black,above right] at ( .60,.36,.216) {$\gamma(t)$};
\end{tikzpicture}
\caption{The box centered at $c_Q=\gamma(0)$ for the moment curve
$\gamma(t)=(t,t^2,t^3)$, with $|t|\le\delta$ and $C=1$. The side lengths
in the directions $e_1,e_2,e_3$ are $2\delta,2\delta^2,2\delta^3$, respectively.}
\label{fig:moment-box}
\end{figure}

For the moment curve, the adapted directions of the boxes are those of the Frenet frame, which Guo, Iosevich, Zhang, and Zorin-Kranich \cite{GIZZK23} use to orient their boxes.
For smooth graphs, Cao, Li, Miao, and Pang instead enclose the convex hull of each graph piece in a rectangular box \cite{CLMP25}.
For illustration, Figure~\ref{fig:moment-box} shows the box centered at $\gamma(0)$, where $\gamma:\R\to\R^3$ is the moment curve, defined by $\gamma(t)=(t,t^2,t^3)$ for $t\in\R$.

To build the cutoff, we first need boxes that enclose the images of the parameter balls.
The next lemma shows that the same enlargement factor works for every center and all sufficiently small scales.
We then construct a smooth function supported in a slightly larger box for each piece and combine these functions into a cutoff equal to one near the Fourier support.

\begin{lemma}\label{lem:taylor}
Fix $C_0>0$.
There exist $0<\delta_0<1$ and $C_1>0$, depending only on $C_0$, $\Omega$, and the parametrization $F$, with the following property.
For each $u_0\in\Omega$, one can choose an orthogonal matrix $A_{u_0}\in O(d)$ whose columns form an adapted basis as above such that if $0<\delta\le\delta_0$ and $u\in\R^m$ satisfies $|u-u_0|\le C_0\delta$, then $u\in U$ and
\begin{equation}\label{eq:taylor-box}
 F(u)\in Q(u_0,\delta;C_1)=F(u_0)+A_{u_0}D_\delta[-C_1,C_1]^d.
\end{equation}
\end{lemma}
\begin{proof}
Here we use compactness to choose $\varepsilon>0$ so that the closed $\varepsilon$-neighborhood $W=\{v\in\R^m:\operatorname{dist}(v,\Omega)\le\varepsilon\}$ is contained in $U$.
With this neighborhood fixed, choose $0<\delta_0<\min\{1,\varepsilon/C_0\}$.
For $v\in W$ and each integer $1\le q\le k$, let $D^qF(v):(\R^m)^q\to\R^d$ be its $q$th Fr\'echet derivative at $v$.
Set
\[
 M_q=\sup_{v\in W}\|D^qF(v)\|_{\mathrm{op}},
 \qquad C_1=1+\max_{1\le q\le k}\frac{M_qC_0^q}{q!}.
\]
The constants $M_q$ are finite by continuity and compactness, so $C_1$ is independent of the center and scale.

Fix $u_0\in\Omega$ and take the adapted orthonormal basis $(e_1,\ldots,e_d)$ from the box construction: its first $r_j$ vectors span $V^{(j)}(F(u_0))$ for each $1\le j\le k$.
These vectors give the side directions of the boxes centered at $F(u_0)$ and form the columns of $A_{u_0}$.
Let $0<\delta\le\delta_0$ and let $u\in\R^m$ such that $|u-u_0|\le C_0\delta$.
Since $\delta<\varepsilon/C_0$, one has $|u-u_0|\leq \varepsilon$, and hence $u\in W\subset U$.
The entire segment from $u_0$ to $u$ also lies in $W$, since every point of that segment is within $\varepsilon$ of $u_0\in\Omega$.

Next, fix an index $1\le\ell\le d$ and put $j=w_\ell$, so $r_{j-1}<\ell\le r_j$.
It remains to show that the displacement along $e_\ell$ is bounded in absolute value by $C_1\delta^j$.
The scalar function
\[
 g_\ell(s)=\la F(u_0+s(u-u_0))-F(u_0),e_\ell\ra,\qquad 0\le s\le1,
\]
records the displacement along this axis, with $g_\ell(0)=0$.
Set $h=u-u_0$.
For every integer $1\le q<j$, the vector $D^qF(u_0)[h,\ldots,h]$ is a linear combination of partial derivatives of order $q$, hence belongs to the $(j-1)$th-order osculating space $V^{(j-1)}(F(u_0))$.
Since $e_\ell$ is orthogonal to that space, the chain rule gives
\[
 g_\ell^{(q)}(0)=\big\langle D^qF(u_0)[h,\ldots,h],e_\ell\big\rangle=0.
\]
Also, $|g_\ell^{(j)}(s)|\le M_j|h|^j$ on $[0,1]$, since the segment lies in $W$ and $|e_\ell|=1$.
Taylor's formula therefore gives
\[
 |g_\ell(1)|
 =\left|\frac{1}{(j-1)!}\int_0^1(1-s)^{j-1}g_\ell^{(j)}(s)\,ds\right|
 \le \frac{M_j}{j!}|h|^j
 \le C_1\delta^{w_\ell}.
\]
These are the bounds on the coordinates of $A_{u_0}^{T}(F(u)-F(u_0))$, which give \eqref{eq:taylor-box}.
This proves the lemma.
\end{proof}

The boxes from Lemma~\ref{lem:taylor} give the right volume at each scale, but their orientations vary with the center.
To control overlap, we retain the parameter cubes that index the boxes.
The cutoff in frequency space will include a factor supported over its parameter cube enlarged by a fixed constant.

At a fixed scale $\delta>0$, one can use a partition of unity in the parameter variables, as Duan, Pham, and Ryou do for the moment curve \cite[Section~2]{DPR25}.
A Hausdorff-dimension bound, however, gives covers at varying scales and need not control the number of cubes at any particular scale.
Bounded overlap at each scale does not give bounded overlap across scales.
We therefore construct cutoffs whose inverse Fourier transforms have bounded $L^1(\R^d)$ norms, uniformly in the center and scale of the box.

Set $c=3/2$ and choose $C_0>c\sqrt m/2$.
Take $\delta_0$ and $C_1$ from Lemma~\ref{lem:taylor} for this choice of $C_0$.
Let $I\subset\R^m$ be a cube centered at $u_I\in\Omega$ with side length $\ell(I)=\delta$, where $0<\delta\le\delta_0$.
For convenience, write $cI$ for the concentric cube of side length $c\delta$, so that $cI\subset B(u_I,C_0\delta)$.
Let $\psi_I\in C_c^\infty(\R^m)$ be supported in $cI$ and satisfy
\begin{equation}\label{eq:psi-bounds}
 |\partial^\nu\psi_I|\le C_\nu\delta^{-|\nu|}
 \quad\text{for every }\nu\in\N_0^m,
\end{equation}
where $C_\nu$ is independent of $I$ and $\delta$.

Since $\widehat f$ is a distribution on the ambient frequency space $\R^d$, we define the cutoff on that space rather than on the parameter domain.
We write $\xi=(\xi',\xi'')\in\R^m\times\R^{d-m}$, so that $\pi(\xi)=\xi'$.
To localize in frequency space, choose $\rho\in C_c^\infty(\R^d)$ with $0\le\rho\le1$ and $\rho=1$ on a neighborhood of $[-C_1,C_1]^d$.
Write $Q_I=Q(u_I,\delta;C_1)$ for the anisotropic box centered at $F(u_I)$.
By Lemma~\ref{lem:taylor}, $F(cI)\subset Q_I$.
The rescaled $\rho$ is one on a neighborhood of $Q_I$.
We define the cutoff in frequency space associated with $Q_I$ by
\begin{equation}\label{eq:bI}
 b_{Q_I}(\xi)=\psi_I(\xi')\,\rho\bigl(D_\delta^{-1}A_{u_I}^{T}(\xi-F(u_I))\bigr) \quad \xi\in \R^d.
\end{equation}
The factor $\psi_I(\xi')$ makes $b_{Q_I}(\xi)$ vanish unless $\xi'\in\supp\psi_I$, so we can count overlapping cutoffs by counting the parameter supports that contain $\xi'$.
The rescaled $\rho$ vanishes outside a fixed enlargement of $Q_I$.

\begin{lemma}[Uniform estimates for a single cutoff]\label{lem:symbol}
Let $I\subset\R^m$ be a cube centered at $u_I\in\Omega$ with side length $0<\delta\le\delta_0$.
Let $\psi_I\in C_c^\infty(\R^m)$ be supported in $cI$ and satisfy \eqref{eq:psi-bounds}, and define $b_{Q_I}$ by \eqref{eq:bI}.
Then
\begin{align}
 \norm{\check b_{Q_I}}_1&\lesssim1,
 &\norm{b_{Q_I}}_\infty&\lesssim1,
 &|\supp b_{Q_I}|&\lesssim\delta^{P_{\mathcal M}}.\label{eq:single-estimates}
\end{align}
The implicit constants are independent of the center $u_I$ and the scale $\delta$.
\end{lemma}
\begin{proof}
Fix the cube $I$ of side length $\delta$ and its cutoff $b_{Q_I}$ from \eqref{eq:bI}.
To estimate $b_{Q_I}$ uniformly in $u_I$ and $\delta$, we move $Q_I$ to the origin and rescale along its sides.
More precisely, set $A_I=A_{u_I}D_\delta$ and write $\xi=F(u_I)+A_Iz$, where $z\in\R^d$.

Since $ \pi(F(u_I))=u_I$, we may write the normalized cutoff as
\[
 B_I(z):=b_{Q_I}(F(u_I)+A_Iz)=\psi_I(u_I+\pi A_Iz)\rho(z),\qquad z\in\R^d.
\]
The support of $B_I$ lies in the fixed compact set $\supp\rho$.
The change of variables also makes the derivative bounds uniform: each derivative falling on $\psi_I(u_I+\pi A_Iz)$ contributes a factor of $\pi A_I$, whose operator norm is at most $\delta$.
Thus the chain rule and \eqref{eq:psi-bounds} give, for every multi-index $\gamma\in\N_0^d$,
\[
 \bigl\|\partial_z^\gamma[\psi_I(u_I+\pi A_Iz)]\bigr\|_\infty
 \lesssim_\gamma \delta^{|\gamma|}\delta^{-|\gamma|}=1.
\]
By the product rule, the derivatives of $B_I$ are uniformly bounded at every fixed order.
In particular, $\norm{b_{Q_I}}_\infty=\norm{B_I}_\infty\lesssim1$, and $|\supp b_{Q_I}|\le |\det A_I|\,|\supp\rho|\lesssim\delta^{P_{\mathcal M}}$.

Next, fix an integer $N>d/2$ and write $\Delta_z$ for the Laplacian in the normalized frequency variable $z$.
To obtain decay in the dual variable $y\in\R^d$, integrate by parts:
\[
 (1+4\pi^2|y|^2)^N|\check B_I(y)|
 \le \norm{(1-\Delta_z)^N B_I}_1\lesssim_N1,
 \qquad y\in\R^d.
\]
Since $2N>d$, the decay bound gives $\|\check B_I\|_1\lesssim 1$, uniformly in $I$.
To transfer the $L^1$ bound for $\check B_I$ to $\check b_{Q_I}$, substitute $\xi=F(u_I)+A_Iz$ in the inverse Fourier integral.
For $x\in\R^d$, we have $d\xi=|\det A_I|\,dz$ and $x\cdot\xi=x\cdot F(u_I)+(A_I^Tx)\cdot z$.
Hence
\begin{align*}
 \check b_{Q_I}(x)
 &=\ee{x\cdot F(u_I)}|\det A_I|\int_{\R^d}B_I(z)\ee{(A_I^Tx)\cdot z}\,dz\\
 &=\ee{x\cdot F(u_I)}|\det A_I|\check B_I(A_I^Tx).
\end{align*}
The exponential has modulus one.
In the integral of the absolute value, set $y=A_I^Tx$, so that $dy=|\det A_I|\,dx$.
Thus
\[
 \norm{\check b_{Q_I}}_1
 =\int_{\R^d}|\det A_I|\,|\check B_I(A_I^Tx)|\,dx
 =\int_{\R^d}|\check B_I(y)|\,dy\lesssim1.
\]
This completes \eqref{eq:single-estimates}.
This completes the proof of the lemma.
\end{proof}
The next lemma bounds the sum of the cutoff kernels when the cubes have comparable side lengths.
\begin{lemma}[Cutoff estimates at one scale]\label{lem:assembly}
Fix $0<\delta\le\delta_0$ and let $\mathcal I$ be a family of $n$ cubes in $\R^m$ centered in $\Omega$.
Assume $0<\ell(I)\le\delta_0$ and $\ell(I)\approx\delta$ for every $I\in\mathcal I$, with fixed comparison constants.
For each $I\in\mathcal I$, take the cutoff $\psi_I$ in parameter space and the cutoff $b_{Q_I}$ in frequency space from \eqref{eq:psi-bounds}--\eqref{eq:bI} at scale $\ell(I)$.
Suppose every $u\in\R^m$ lies in at most $M_0$ of the supports $\supp\psi_I$, where $M_0$ is fixed.
Set $\chi=\sum_{I\in\mathcal I}b_{Q_I}$ and $K=\check\chi$.
Then
\begin{align}
 \norm K_1&\lesssim n,&
 \norm K_2&\lesssim(n\delta^{P_{\mathcal M}})^{1/2},&
 \norm K_\infty&\lesssim n\delta^{P_{\mathcal M}}.\label{eq:assembly-basic}
\end{align}
Moreover, if $2<p<\infty$ and $p'=p/(p-1)$, then
\begin{equation}\label{eq:assembly-lorentz}
 \norm K_{p',1}\lesssim_p n^{1-1/p}\delta^{P_{\mathcal M}/p}.
\end{equation}
\end{lemma}
\begin{proof}
Recall that $K=\sum_{I\in\mathcal I}\check b_{Q_I}$ and that the cubes have side lengths comparable to $\delta$.
By the triangle inequality and \eqref{eq:single-estimates}, $\norm K_1\le\sum_{I\in\mathcal I}\norm{\check b_{Q_I}}_1\lesssim n$.

For the $L^2$ bound, by contrast, fix $\xi=(\xi',\xi'')\in\R^m\times\R^{d-m}$.
Since $b_{Q_I}(\xi)$ vanishes unless $\xi'\in\supp\psi_I$, at most $M_0$ summands of $\chi(\xi)$ are nonzero.
Cauchy--Schwarz therefore gives
\[
 |\chi(\xi)|^2\le M_0\sum_{I\in\mathcal I}|b_{Q_I}(\xi)|^2.
\]
Integrate this inequality and apply Plancherel.
The support and supremum bounds in \eqref{eq:single-estimates} give $\norm K_2^2\lesssim M_0 n\delta^{P_{\mathcal M}}$.
Moreover, the same bounds give $\norm K_\infty\le\norm\chi_1\lesssim n\delta^{P_{\mathcal M}}$ by Fourier inversion.

Finally, Lorentz interpolation (Lemma~\ref{lem:interpolation}) gives $\norm K_{p',1}\lesssim_p\norm K_1^{1-2/p}\norm K_2^{2/p}$.
Substituting the two bounds proves \eqref{eq:assembly-lorentz}.
This completes the proof of the lemma.
\end{proof}

To apply Lemma~\ref{lem:cutoff}, we further need the sum of the cutoffs in frequency space to equal one on an open neighborhood of $E$, not just on $E$ itself.
This ensures that multiplication by the sum leaves the distribution $\widehat f$ unchanged.
The next lemma obtains such a neighborhood from a partition of unity in parameter space, even when the boxes have different scales.

\begin{lemma}[Cutoffs equal to one near the graph]\label{lem:cutoff-neighborhood}
Let $\mathcal I$ be a finite nonempty family of cubes in $\R^m$ such that $I\in\mathcal I$ is centered at $u_I\in\Omega$ and has side length $0<\ell(I)\le\delta_0$; the side lengths may differ.
For each $I\in\mathcal I$, let $\psi_I$ be a cutoff in parameter space satisfying \eqref{eq:psi-bounds} at scale $\ell(I)$, with $\supp\psi_I\subset B(u_I,C_0\ell(I))$, where $C_0$ and $\delta_0$ are as in Lemma~\ref{lem:taylor}.
Let $b_{Q_I}$ be the corresponding cutoffs in frequency space defined by \eqref{eq:bI}.
If $\sum_{I\in\mathcal I}\psi_I=1$ on an open set $ W\subset U$, then $\sum_{I\in\mathcal I} b_{Q_I}=1$ on an open neighborhood of $ F(W)$ in $\R^d$.
\end{lemma}

\begin{proof}
Recall that $F(u)=(u,\Phi(u))$ for $u\in U$.
For every $I\in\mathcal I$, recall from \eqref{eq:bI} that
\[
 b_{Q_I}(\xi)=\psi_I(\xi')\,\rho\bigl(D_{\ell(I)}^{-1}A_{u_I}^{T}(\xi-F(u_I))\bigr),
 \qquad \xi=(\xi',\xi'')\in\R^m\times\R^{d-m}.
\]
We show that there is an open neighborhood of $ F(W)$ such that, for every $\xi$ in this neighborhood and every $I\in\mathcal I$, we have $b_{Q_I}(\xi)=\psi_I(\xi')$.
Here compactness gives a uniform neighborhood: the open set on which $\rho$ equals one contains the $\varepsilon$-neighborhood of $[-C_1,C_1]^d$ for some $\varepsilon>0$.
Set $\eta=\tfrac \varepsilon2\min_{I\in\mathcal I}\ell(I)^k>0$, so that the same vertical distance $\eta$ works for every cube.

To verify the equality, fix $\xi=(\xi',\xi'')\in\R^m\times\R^{d-m}$ with $\xi'\in W$ and $\abs{\xi''-\Phi(\xi')}<\eta$.
If $\psi_I(\xi')=0$, then $b_{Q_I}(\xi)=0$.
Otherwise, set $A_I=A_{u_I}D_{\ell(I)}$.
Since $\xi'\in\supp\psi_I\subset B(u_I,C_0\ell(I))$, Lemma~\ref{lem:taylor} gives
\[
 A_I^{-1}\bigl(F(\xi')-F(u_I)\bigr)\in[-C_1,C_1]^d.
\]
Note that the smallest singular value of $A_I$ is $\ell(I)^k$.
As $\xi-F(\xi')=(0,\xi''-\Phi(\xi'))$, we have
\[
 \abs{A_I^{-1}\bigl(\xi-F(\xi')\bigr)}
 \le \ell(I)^{-k}\abs{\xi''-\Phi(\xi')}<\ell(I)^{-k}\eta<\varepsilon,
\]
by the choice of $\eta$.
Thus $A_I^{-1}(\xi-F(u_I))$ lies in the $\varepsilon$-neighborhood of $[-C_1,C_1]^d$, so the $\rho$ factor in \eqref{eq:bI} is one.
Summing $b_{Q_I}(\xi)=\psi_I(\xi')$ over $I\in\mathcal I$ gives $\sum_I b_{Q_I}(\xi)=1$.
Finally, the conditions $\xi'\in W$ and $\abs{\xi''-\Phi(\xi')}<\eta$ define an open neighborhood of $ F(W)$.
This proves the lemma.
\end{proof}

\subsection{Dyadic covers and cutoff estimates across scales}\label{sec:dyadic-cutoffs}

Hausdorff measure gives a dyadic cover in which the side lengths may vary.
We construct cutoffs on these cubes whose sum is at most one across all scales.
This bound, rather than bounded overlap at each scale, is what we use in the $L^2$ estimate.
The ordered dyadic partition and its scale-uniform derivative estimates are due to Harvey and Polking \cite[Lemma~3.1]{HarveyPolking70}.
We recall the construction in the form needed for the lift to ambient boxes.

\begin{lemma}[Dyadic parameter partition]\label{lem:dyadic-partition}
Let $V\subset\R^m$ be nonempty and compact, and let $W\subset\R^m$ be an open neighborhood of $V$.
Assume that $\HH^s(V)=0$ for some $s>0$.
Let $\varepsilon>0$ and take $\delta_*>0$ sufficiently small in terms of $V$ and $W$.
Then there is a finite dyadic cover $\mathcal I$ of $V$ by cubes with disjoint interiors and $\sum_{I\in\mathcal I}\ell(I)^s<\varepsilon$, such that the following hold.

For every $I\in\mathcal I$, the side length satisfies $\ell(I)\le\delta_*$ and $\overline{(3/2)I}$ is a compact subset of $W$.
For each $I\in\mathcal I$, there is $\psi_I\in C_c^\infty(\R^m)$ such that
\[
\begin{gathered}
 0\leq\psi_I,\qquad \sum_I\psi_I\leq1\ \text{on }\R^m,\qquad
 \sum_I\psi_I=1\ \text{near }V,\\
 \supp\psi_I\subset (3/2)I,\qquad
 |\partial^\nu\psi_I|\le C_{m,\nu}\ell(I)^{-|\nu|}\quad\text{for every }\nu\in\N_0^m.
\end{gathered}
\]

\end{lemma}
\begin{proof}
Fix $s>0$ with $\HH^s(V)=0$ and let $\varepsilon>0$ be the desired upper bound for $\sum_{I\in\mathcal I}\ell(I)^s$.
Take $\delta_*>0$ sufficiently small in terms of $V$ and $W$ so that the enlarged cubes stay in $W$.
This is the upper bound on their original side lengths.
Since $\HH^s(V)=0$, choose a countable dyadic cover $\mathcal J$ of $V$ whose cubes meet $V$, with $\ell(J)\le\delta_*$ for every $J\in\mathcal J$ and $\sum_{J\in\mathcal J}\ell(J)^s<3^{-m}\varepsilon$.

By compactness, choose a finite subfamily $\mathcal J_0\subset\mathcal J$ such that $V\subset\bigcup_{J\in\mathcal J_0}\operatorname{int}(3J)$.
For each $J\in\mathcal J_0$, let $\mathcal D_J$ consist of the $3^m$ closed dyadic cubes of side length $\ell(J)$ that tile $3J$.
Set
\[
 \mathcal D=\bigcup_{J\in\mathcal J_0}\mathcal D_J,
 \qquad
 \mathcal I=\{I\in\mathcal D:\text{ there is no } I'\in\mathcal D\text{ with }I\subsetneq I'\}.
\]
Since every cube in the finite family $\mathcal D$ is contained in one from $\mathcal I$, $\bigcup_{I\in\mathcal I}I=\bigcup_{J\in\mathcal J_0}3J$.
Distinct cubes in $\mathcal I$ have disjoint interiors, since two dyadic cubes with intersecting interiors are nested.
In particular, the retained cubes cover a neighborhood of $V$, and
\[
 \sum_{I\in\mathcal I}\ell(I)^s
 \le\sum_{J\in\mathcal J_0}\sum_{I\in\mathcal D_J}\ell(I)^s
 =3^m\sum_{J\in\mathcal J_0}\ell(J)^s<\varepsilon.
\]
Moreover, each retained cube comes from a $3J$ whose original cube $J$ meets $V$, so $\overline{(3/2)I}\subset W$ if $\delta_*$ is sufficiently small.

Enumerate the retained cubes as $I_1,\ldots,I_N$ so that $\delta_1\ge\cdots\ge\delta_N$, where $\delta_i=\ell(I_i)$.
For each $1\le i\le N$, let $u_{I_i}$ be the center of $I_i$.
Fix $\phi\in C_c^\infty((-3/4,3/4)^m)$ with $0\le\phi\le1$ and $\phi=1$ on $[-1/2,1/2]^m$.
For every $1\le i\le N$, put
\[
 \phi_i(u)=\phi\bigl((u-u_{I_i})/\delta_i\bigr),\qquad
 \psi_{I_i}(u)=\phi_i(u)\prod_{j<i}(1-\phi_j(u)),\qquad u\in\R^m,
\]
with the empty product equal to one; see \cite[proof of Lemma~3.1, pp.~43--44]{HarveyPolking70} for the ordered product construction.
Then $0\le\psi_{I_i}\le1$ and $\supp\psi_{I_i}\subset(3/2)I_i$.
For every $u\in\mathbb R^m$, the sum telescopes:
\[
 \sum_{i=1}^N\psi_{I_i}(u)=1-\prod_{i=1}^N(1-\phi_i(u))
\]
The identity therefore gives a sum bounded by one everywhere and equal to one on $\bigcup_{i=1}^N I_i$, hence near $V$.

To bound derivatives of $\psi_{I_i}$, fix $u\in\R^m$ and $1\le i\le N$.
For a nonzero multi-index $\beta\in\N_0^m$, differentiation gives
\begin{equation}\label{eq:dyadic-bump-derivative}
 |\partial^\beta\phi_j(u)|\le C_{m,\beta}\delta_j^{-|\beta|}\one_{(3/2)I_j}(u).
\end{equation}
Leibniz's rule may place a derivative on any factor $1-\phi_j$ with $j<i$, so we need a bound for $\sum_{j=1}^i|\partial^\beta\phi_j(u)|$ that is independent of $i$.
By \eqref{eq:dyadic-bump-derivative}, only cubes with $u\in(3/2)I_j$ can contribute.
We count these cubes separately at each dyadic scale.

For each integer $h\ge0$, define
\[
 \mathcal J_h(u)=\bigl\{j\in\{1,\ldots,i\}:\delta_j=2^h\delta_i,\ u\in(3/2)I_j\bigr\}.
\]
To count these indices, write the center of each cube of side length $\delta_j=2^h\delta_i$ as $u_{I_j}=2^h\delta_i(n_j+\tfrac12\mathbf 1)$, where $n_j\in\mathbb Z^m$ and $\mathbf 1=(1,\ldots,1)$.
Write $(n_j)_\ell$ and $(u)_\ell$ for the $\ell$th coordinates of $n_j$ and $u$, respectively.
If $j\in\mathcal J_h(u)$, then
\[
 \left|(n_j)_\ell-\left(\frac{(u)_\ell}{2^h\delta_i}-\frac12\right)\right|\le\frac34
 \qquad\text{for every }1\le\ell\le m.
\]
For each coordinate there are at most two integers in this interval of length $3/2$.
Distinct cubes of side length $2^h\delta_i$ have distinct centers, so $|\mathcal J_h(u)|\le 2^m$.

If $1\le j\le i$ and $\partial^\beta\phi_j(u)\ne0$, then $u\in(3/2)I_j$ by the support of $\phi_j$.
Also, $\delta_j/\delta_i$ is a power of two and is at least one, so $\delta_j=2^h\delta_i$ for a unique integer $h\ge0$.
Thus $j\in\mathcal J_h(u)$ for that $h$, and we may group the sum by scale.
Applying \eqref{eq:dyadic-bump-derivative} to each term gives
\begin{equation}\label{eq:dyadic-bump-derivative-sum}
\begin{aligned}
 \sum_{j=1}^i|\partial^\beta\phi_j(u)|
 &=\sum_{h\ge0}\sum_{j\in\mathcal J_h(u)}|\partial^\beta\phi_j(u)|\le C_{m,\beta}\sum_{h\ge0}|\mathcal J_h(u)|(2^h\delta_i)^{-|\beta|}
 \lesssim_{m,\beta}\delta_i^{-|\beta|}.
\end{aligned}
\end{equation}
The geometric series converges because $|\beta|>0$.

We now apply Leibniz's rule to $\psi_{I_i}$.
Fix $\nu\in\N_0^m$ with $|\nu|\ge1$.
In each term, let $1\le j_1<\cdots<j_q\le i$ be the indices of the differentiated factors, and let $\beta_r\in\N_0^m\setminus\{0\}$ be the multi-index applied to the factor with index $j_r$.
Then $1\le q\le|\nu|$ and $\beta_1+\cdots+\beta_q=\nu$.
The undifferentiated factors have absolute value at most one, and a differentiated factor has absolute value $|\partial^{\beta_r}\phi_{j_r}(u)|$, whether it comes from $\phi_i$ or $1-\phi_{j_r}$.
For each fixed decomposition of $\nu$, dropping the ordering and distinctness restrictions gives
\[
 \sum_{1\le j_1<\cdots<j_q\le i}
 \prod_{r=1}^q|\partial^{\beta_r}\phi_{j_r}(u)|
 \le\prod_{r=1}^q\left(\sum_{j=1}^i|\partial^{\beta_r}\phi_j(u)|\right)
 \lesssim_{m,\nu}\delta_i^{-|\nu|},
\]
by \eqref{eq:dyadic-bump-derivative-sum} and $\sum_{r=1}^q|\beta_r|=|\nu|$.
Since the number of decompositions and the Leibniz coefficients depend only on $\nu$, this proves $|\partial^\nu\psi_{I_i}(u)|\lesssim_{m,\nu}\delta_i^{-|\nu|}$, uniformly in $i$, $u$, and the cover.
The case $\nu=0$ follows from $0\le\psi_{I_i}\le1$.
This completes the proof of the lemma.
\end{proof}
\begin{remark}
    \label{rem:finite-dyadic-cost}More generally, the same proof applies when $\HH^s(V)<\infty$. Start with a dyadic cover of cost at most $C_{m,s}(\HH^s(V)+\varepsilon)$ to obtain $\sum_{I\in\mathcal I}\ell(I)^s\le C'_{m,s}(\HH^s(V)+\varepsilon)$, uniformly in $\delta_*$.
\end{remark}
The dyadic cubes $I_i$ from Lemma~\ref{lem:dyadic-partition} may have different side lengths.
Accordingly, the next lemma estimates the sum of the cutoff kernels in terms of the volumes of the individual anisotropic boxes $Q_{I_i}$ constructed using Lemma~\ref{lem:taylor}.
For a related argument with isotropic cutoffs, see \cite[Lemma~4.1]{Dob26}.

\begin{lemma}[Assembly across scales]\label{lem:all-scale-assembly}
Let $I_1,\ldots,I_N$ and $\psi_{I_1},\ldots,\psi_{I_N}$ be the cubes and cutoffs in parameter space from Lemma~\ref{lem:dyadic-partition} or Remark~\ref{rem:finite-dyadic-cost}.
Assume the cube centers lie in a fixed compact set $\Omega\subset U$ and that their side lengths are at most the corresponding $\delta_0$ in Lemma~\ref{lem:taylor}.
For each $1\le i\le N$, let $\delta_i=\ell(I_i)$ and let $Q_{I_i}=Q(u_{I_i},\delta_i;C_1)$ be the anisotropic box from Lemma~\ref{lem:taylor}, centered at $F(u_{I_i})$ and adapted to the osculating spaces of $\mathcal M$ there.
Write $b_i=b_{Q_{I_i}}$ for the cutoff in frequency space associated with this box, defined by \eqref{eq:bI} using $\psi_{I_i}$.
If $1<r<2$ and $c_1,\ldots,c_N\in\C$, then
\begin{equation}\label{eq:all-scale-assembly-lr}
 \Big\|\sum_i c_i\check b_i\Big\|_r^r
 \lesssim_r \sum_i |c_i|^r \delta_i^{P_{\mathcal M}(r-1)},\quad\quad \big\|\sum_i\check b_i\big\|_\infty\lesssim\sum_i \delta_i^{P_{\mathcal M}}.
\end{equation}

\end{lemma}

\begin{proof}
Recall that $b_i=b_{Q_{I_i}}$ is the cutoff associated with the box $Q_{I_i}$ and that $\delta_i=\ell(I_i)$.
To retain the factor $\delta_i^{P_{\mathcal M}(r-1)}$ for each cube $I_i$, we place the $L^1$ and $L^2$ estimates on the same weighted coefficient space.
Let $\nu$ be the measure on $\{1,\ldots,N\}$ with $\nu(\{i\})=\delta_i^{-P_{\mathcal M}}$, and define $T:\C^N\to\Sch(\R^d)$ by
\[
 T(\widetilde{\mathbf c})=\sum_{i=1}^N\frac{\widetilde c_i}{\delta_i^{P_{\mathcal M}}}\check b_i,
 \qquad \widetilde{\mathbf c}=(\widetilde c_1,\ldots,\widetilde c_N)\in\C^N.
\]
First, the uniform $L^1$ bound for $\check b_i$ in \eqref{eq:single-estimates} gives
\[
 \norm{T(\widetilde{\mathbf c})}_1\le\sum_{i=1}^N\frac{|\widetilde c_i|}{\delta_i^{P_{\mathcal M}}}\norm{\check b_i}_1
 \lesssim\norm{\widetilde{\mathbf c}}_{L^1(\nu)}.
\]
For the $L^2$ bound, the key point is that the cutoffs sum to at most one even though the boxes have different sizes.
Indeed, writing $\xi=(\xi',\xi'')\in\R^m\times\R^{d-m}$, we have $0\le b_i(\xi)\le\psi_{I_i}(\xi')$ by \eqref{eq:bI}.
Thus $\sum_{i=1}^N b_i(\xi)\le1$ by Lemma~\ref{lem:dyadic-partition}.
By Plancherel and weighted Cauchy--Schwarz,
\begin{align*}
 \norm{T(\widetilde{\mathbf c})}_2^2
 &\le\int_{\R^d}\Big(\sum_{i=1}^N b_i(\xi)\Big)
          \sum_{i=1}^N\frac{|\widetilde c_i|^2}{\delta_i^{2P_{\mathcal M}}}b_i(\xi)\dd\xi\\
 &\le\sum_{i=1}^N\frac{|\widetilde c_i|^2}{\delta_i^{2P_{\mathcal M}}}\int_{\R^d}b_i(\xi)\dd\xi
 \lesssim\sum_{i=1}^N\frac{|\widetilde c_i|^2}{\delta_i^{P_{\mathcal M}}}
 =\norm{\widetilde{\mathbf c}}_{L^2(\nu)}^2.
\end{align*}
Here the last inequality uses $\int_{\R^d}b_i(\xi)\dd\xi\lesssim\delta_i^{P_{\mathcal M}}$, which follows from the support and supremum bounds in \eqref{eq:single-estimates}.
Interpolating these two estimates by Riesz--Thorin \cite[Theorem~1.3.4]{Grafakos14}, we obtain $\norm{T(\widetilde{\mathbf c})}_r^r\lesssim_r\norm{\widetilde{\mathbf c}}_{L^r(\nu)}^r$ for $1<r<2$.
Taking $\widetilde c_i=\delta_i^{P_{\mathcal M}}c_i$ gives $T(\widetilde{\mathbf c})=\sum_{i=1}^N c_i\check b_i$ and
\[
 \norm{\widetilde{\mathbf c}}_{L^r(\nu)}^r=\sum_{i=1}^N |c_i|^r\delta_i^{P_{\mathcal M}(r-1)}.
\]
This proves the first estimate.
For the $L^\infty$ estimate, Fourier inversion and the same integral bound give
\[
 \Big\|\sum_{i=1}^N\check b_i\Big\|_\infty
 \le\sum_{i=1}^N\int_{\R^d}b_i(\xi)\dd\xi
 \lesssim\sum_{i=1}^N\delta_i^{P_{\mathcal M}}.
\]
This completes the proof of the lemma.
\end{proof}

\subsection{Proofs of the Hausdorff dimension and finite-measure theorems}\label{sec:hausdorff-proofs}

\begin{proof}[Proof of Theorem~\ref{thm:hausdorff}]
Let $d,m,k$ be integers with $d\ge2$, $1\le m<d$, and $k\ge2$.
Recall that $\mathcal M=F(U)\subset\R^d$ is an $m$-dimensional $C^k$ graph of type $k$ with common rank profile $(r_1,\ldots,r_k)$ and volume exponent $P_{\mathcal M}=\sum_{j=1}^k j(r_j-r_{j-1})$, where $r_0=0$.
Let $E\subset\mathcal M$ be compact with $\dimH E=\alpha\in(0,m]$, and suppose $f\in L^{p,\infty}(\R^d)$ represents a tempered distribution with $\supp\widehat f\subset E$.
We show that $f=0$ whenever $1\le p<p_{\mathcal M,\alpha}=1+P_{\mathcal M}/\alpha$.

We first choose a strong Lebesgue exponent.
Since $\alpha\le m<d<P_{\mathcal M}$ and $p<p_{\mathcal M,\alpha}$, we can choose $p_0$ with
\[
 \max\{p,1+P_{\mathcal M}/d\}<p_0<p_{\mathcal M,\alpha}.
\]
Set $s=P_{\mathcal M}/(p_0-1)$.
The two bounds on $p_0$ give $\alpha<s<d<P_{\mathcal M}$.
Since $\dimH E=\alpha$, it follows that $\HH^s(E)=0$.
Moreover, $p_0>p$, so $f\in L^{p_0}(\R^d)$ by \eqref{eq:bandlimited-weak-l1} if $p=1$ and by Lemma~\ref{lem:lorentz}(iii) if $p>1$.

Write $\pi:\R^d\to\R^m$ for projection onto the first $m$ coordinates and set $V=F^{-1}(E)=\pi(E)\subset U$.
Since $\pi$ is $1$-Lipschitz, $\HH^s(V)\le\HH^s(E)=0$.
We use this to cover $V$ by dyadic cubes and build cutoffs that equal one on a neighborhood of $E$.
Choose an open set $W\subset\R^m$ containing $V$ with $\overline W$ compact in $U$.

Applying Lemma~\ref{lem:taylor} with $\Omega=\overline W$, there are constants $0<\delta_0<1$ and $C_1>0$ such that \eqref{eq:taylor-box} holds for every $u_0\in\Omega$ and every $0<\delta\le\delta_0$.
The box estimate sets an upper bound on the side lengths, but the cubes need not all have the same size.

Given a sufficiently small $\delta_*\in(0,\delta_0]$, apply Lemma~\ref{lem:dyadic-partition} with $\varepsilon=1$ to obtain a dyadic cover $I_1,\ldots,I_N$ of $V$ in $\R^m$.
Writing $\delta_i=\ell(I_i)$, we have $\delta_i\le\delta_*$ and $(3/2)I_i\subset W\subset\Omega$ for every $1\le i\le N$.
In addition, $\sum_{i=1}^N\delta_i^s<1$.
To turn this cover into a smooth cutoff near $V$, we use the functions $\psi_{I_i}\in C_c^\infty(\R^m)$ supplied by the lemma.
For every $1\le i\le N$ and $\nu\in\N_0^m$,
\begin{equation}\label{eq:hausdorff-cutoff-properties}
 0\le\psi_{I_i},\qquad \supp\psi_{I_i}\subset(3/2)I_i,\qquad
 \norm{\partial^\nu\psi_{I_i}}_\infty\le C_{m,\nu}\delta_i^{-|\nu|}.
\end{equation}
Although the supports can overlap, their sum is at most one on $\R^m$.
Moreover, the sum equals one on an open neighborhood of $V$ contained in $W$.

We now pass from the cutoffs in parameter space to cutoffs in frequency space.
For each $1\le i\le N$, let $u_{I_i}$ be the center of $I_i$ and let $Q_{I_i}=Q(u_{I_i},\delta_i;C_1)\subset\R^d$ be the adapted box centered at $F(u_{I_i})$.
By \eqref{eq:hausdorff-cutoff-properties} and the choice $C_0>3\sqrt m/4$, we have
\begin{equation*}
 \supp\psi_{I_i}\subset(3/2)I_i\subset B(u_{I_i},C_0\delta_i)
 \qquad\text{for every }1\le i\le N.
\end{equation*}
Recall from \eqref{eq:bI} that the cutoff $b_{Q_{I_i}}\in C_c^\infty(\R^d)$ adapted to $Q_{I_i}$ is given by
\[
 b_{Q_{I_i}}(\xi)=\psi_{I_i}(\pi(\xi))\,
 \rho\bigl(D_{\delta_i}^{-1}A_{u_{I_i}}^T(\xi-F(u_{I_i}))\bigr),\qquad \xi\in\R^d.
\]
Since $\sum_{i=1}^N\psi_{I_i}=1$ near $V$, Lemma~\ref{lem:cutoff-neighborhood} gives
\[
 \sum_{i=1}^N b_{Q_{I_i}}(\xi)=1
 \qquad\text{for every $\xi$ in an open neighborhood of $E=F(V)$ in $\R^d$.}
\]
To apply the vanishing criterion in Lemma~\ref{lem:cutoff}(ii), we need a uniform $L^{p_0'}$ bound and a vanishing $L^\infty$ bound for the inverse Fourier transform of the cutoff.
Set $\chi_{\delta_*}=\sum_{i=1}^N b_{Q_{I_i}}$ and $K_{\delta_*}=\check\chi_{\delta_*}$.
Since $\chi_{\delta_*}\in C_c^\infty(\R^d)$, its inverse Fourier transform belongs to $\Sch(\R^d)$.

The choice $p_0=1+P_{\mathcal M}/s$ makes the exponent in the assembly estimate match the covering cost.
Indeed, $p_0'=1+s/P_{\mathcal M}$ lies in $(1,2)$ and satisfies $P_{\mathcal M}(p_0'-1)=s$.
Applying Lemma~\ref{lem:all-scale-assembly} with $r=p_0'$ and $c_i=1$ for every $1\le i\le N$ gives the first bound below.
For the second, use the $L^\infty$ bound in \eqref{eq:all-scale-assembly-lr}  and $\delta_i\le\delta_*$:
\begin{equation}\label{eq:hausdorff-kernel-bounds}
 \norm{K_{\delta_*}}_{p_0'}^{p_0'}\lesssim\sum_{i=1}^N\delta_i^s,
 \qquad
 \norm{K_{\delta_*}}_\infty\lesssim\sum_{i=1}^N\delta_i^{P_{\mathcal M}}
 \le\delta_*^{P_{\mathcal M}-s}\sum_{i=1}^N\delta_i^s.
\end{equation}
The sum $\sum_{i=1}^N\delta_i^s$ is bounded independently of $\delta_*$.
Therefore the first norm stays bounded, whereas the second tends to zero as $\delta_*\downarrow0$, since $P_{\mathcal M}-s>0$.

Choose a sequence of admissible scale bounds $\delta_n\downarrow0$ and set $\chi_n=\chi_{\delta_n}$ and $K_n=K_{\delta_n}$.
Each $\chi_n$ equals one on a neighborhood of $E$, and \eqref{eq:hausdorff-kernel-bounds} gives $\sup_n\norm{K_n}_{p_0'}<\infty$ and $\norm{K_n}_\infty\to0$.
Thus Lemma~\ref{lem:cutoff}(ii) applies to $f\in L^{p_0}(\R^d)$ with $\supp\widehat f\subset E$ and yields $f=0$.

This proves Theorem~\ref{thm:hausdorff} on a graph patch.
\end{proof}

\begin{proof}[Proof of Theorem~\ref{thm:finite-hausdorff}]
Let $E\subset\mathcal M$ be compact with $\HH^\alpha(E)<\infty$, where $0<\alpha\le m$, and suppose $f\in L^p(\R^d)$ with $\supp\widehat f\subset E$.
Set $p_0=p_{\mathcal M,\alpha}$.
It suffices to prove that $f=0$ when $f\in L^{p_0}(\R^d)$.
Indeed, if $p<p_0$, compact Fourier support gives $f\in L^{p_0}(\R^d)$ by \eqref{eq:bandlimited-weak-l1} for $p=1$ and by Lemma~\ref{lem:lorentz}(iii) for $p>1$, using $L^p\subset L^{p,\infty}$.

As in the preceding proof, we construct a smooth cutoff equal to one near $E$ and estimate its inverse Fourier transform.
The difference is that at $p_0=p_{\mathcal M,\alpha}$ the covering exponent is $s=\alpha$, so we know that the relevant Hausdorff measure is finite, not necessarily zero.
Since $V=F^{-1}(E)=\pi(E)$ and $\pi$ is $1$-Lipschitz, $\HH^\alpha(V)\le\HH^\alpha(E)<\infty$.
Hence Remark~\ref{rem:finite-dyadic-cost} gives, for every sufficiently small $\delta_*>0$, dyadic cubes $I_1,\ldots,I_N$ covering $V$ with side lengths $\delta_i=\ell(I_i)\le\delta_*$ and
\[
 \sum_{i=1}^N\delta_i^\alpha\le C_{m,\alpha}\bigl(1+\HH^\alpha(V)\bigr).
\]
The associated cutoffs have the same properties as in the preceding proof.
In particular, their sum in frequency space, $\chi_{\delta_*}$, equals one on a neighborhood of $E$.

The key point is that the covering cost only needs to be bounded uniformly in $\delta_*$.
Since $p_0'=1+\alpha/P_{\mathcal M}$, \eqref{eq:hausdorff-kernel-bounds} with $s=\alpha$ gives, for $K_{\delta_*}=\check\chi_{\delta_*}$,
\[
 \norm{K_{\delta_*}}_{p_0'}^{p_0'}\lesssim 1+\HH^\alpha(V),
 \qquad
 \norm{K_{\delta_*}}_\infty\lesssim\delta_*^{P_{\mathcal M}-\alpha}\bigl(1+\HH^\alpha(V)\bigr).
\]
The first bound is uniform, while the second tends to zero as $\delta_*\downarrow0$ because $\alpha\le m<P_{\mathcal M}$.
Thus, along any sequence of admissible $\delta_*\downarrow0$, Lemma~\ref{lem:cutoff}(ii) gives $f=0$.

For later use, the endpoint argument also applies with $s$ in place of $\alpha$: the only restriction on $s$ in the cutoff estimates is $0<s<P_{\mathcal M}$.
In particular,
\begin{equation}\label{eq:finite-hausdorff-uniqueness}
 \HH^s(E)<\infty,\quad g\in L^{1+P_{\mathcal M}/s}(\R^d),\quad
 \supp\widehat g\subset E \quad\Longrightarrow\quad g=0,
 \qquad 0<s<P_{\mathcal M}.
\end{equation}
This proves the endpoint criterion on a graph patch.
\end{proof}

We now apply Section~\ref{sec:chart-localization} to pass from a graph patch to $\mathcal M$.
If $g$ satisfies the hypotheses of \eqref{eq:finite-hausdorff-uniqueness}, then Young's inequality gives
\[
 \norm{g*\check\eta_j}_{p_0}\le\norm{g}_{p_0}\norm{\check\eta_j}_1,
 \qquad \HH^s(E_j)\le\HH^s(E),
 \qquad p_0=1+P_{\mathcal M}/s.
\]
Thus \eqref{eq:finite-hausdorff-uniqueness} applies to each localized term.
Summing gives $g=0$, so the criterion holds on $\mathcal M$.
Taking $s=\alpha$ gives the strong endpoint.
In the strict weak range, choose $p_0$ and $s$ as in the proof of Theorem~\ref{thm:hausdorff}.
The same bandlimited upgrade gives $f\in L^{p_0}(\R^d)$ before localization, and $\HH^s(E)=0$.
The global criterion therefore gives $f=0$.
This completes the proofs of Theorems~\ref{thm:hausdorff} and~\ref{thm:finite-hausdorff}.

\section{Proof of Theorem~\ref{thm:covering}}\label{sec:proof-covering}

In this section we prove Theorem~\ref{thm:covering}.
The proof of Theorem~\ref{thm:hausdorff} combines cutoffs at different scales.
Here, by contrast, we use a cover by balls of a common radius and construct all the cutoffs at that scale.
This fixed-scale construction follows the cutoff and interpolation argument in \cite[proof of Proposition~7.13]{CLMP25}.
Below the critical exponent, the Lorentz norm of the resulting kernel has an extra positive power of the radius, so bounded normalized covering numbers along a sequence suffice.
At the endpoint,  however, this power disappears, and we instead use a sequence along which the normalized covering number tends to zero.

\begin{proof}[Proof of Theorem~\ref{thm:covering}]
Let $E\subset\mathcal M$ be compact and fix $0<\alpha\le m$.
Suppose $f\in L^{p,\infty}(\R^d)$, where $1\le p<\infty$, represents a tempered distribution with $\supp\widehat f\subset E$.
We will prove that $f=0$ below $p_{\mathcal M,\alpha}$ under \eqref{eq:lower-cover-assumption}, and at $p=p_{\mathcal M,\alpha}$ under \eqref{eq:zero-content}.

For $p>2$, Section~\ref{sec:chart-localization} and Lemma~\ref{lem:lorentz}(ii) give $f_j\in L^{p,\infty}(\R^d)$.
Since
\[
 \delta^\alpha\cov(E_j,\delta)\le\delta^\alpha\cov(E,\delta)
 \qquad (\delta>0),
\]
the same sequence of scales works on every patch.
Thus for $p>2$ we may assume that $\mathcal M=F(U)$ is a graph patch, with the same rank profile $(r_1,\ldots,r_k)$ and volume exponent $P_{\mathcal M}$.
For $1\le p\le2$, we use the global criterion \eqref{eq:finite-hausdorff-uniqueness} below.

If $E$ is empty, then $\widehat f=0$, so assume $E$ is nonempty.

\smallskip\noindent\textit{Part (i).} Assume $1\le p<p_{\mathcal M,\alpha}$ and recall the covering hypothesis \eqref{eq:lower-cover-assumption}:
\begin{equation}\label{eq:covering-proof-bounded-content}
 \liminf_{\delta\downarrow0}\delta^\alpha\cov(E,\delta)<\infty.
\end{equation}
For $1\le p\le2$, no covering-number assumption is needed.
Since $m<d<P_{\mathcal M}$, choose $1+P_{\mathcal M}/d<p_0<p_{\mathcal M,m}$ and set $s=P_{\mathcal M}/(p_0-1)$.
Then $p_0>2$ and $m<s<d<P_{\mathcal M}$, so $\HH^s(E)=0$ because $\dimH E\le m$.
Moreover, compact Fourier support gives $f\in L^{p_0}(\R^d)$ by \eqref{eq:bandlimited-weak-l1} if $p=1$ and by Lemma~\ref{lem:lorentz}(iii) if $p>1$.
Thus \eqref{eq:finite-hausdorff-uniqueness} gives $f=0$.

Suppose now that $2<p<p_{\mathcal M,\alpha}$ and write $p'=p/(p-1)$.
For each $\delta>0$, set $A_\delta=\delta^\alpha\cov(E,\delta)$.
By \eqref{eq:covering-proof-bounded-content}, choose a sequence $\delta_n\downarrow0$ along which $A_{\delta_n}$ stays bounded.
To apply Lemma~\ref{lem:cutoff}(i), we will construct, for every sufficiently small $\delta$, a frequency cutoff equal to one near $E$ and estimate its inverse transform in $L^{p',1}(\R^d)$.
We then use this estimate at the chosen scales.
The construction itself does not use \eqref{eq:lower-cover-assumption}, so it will also apply in part~(ii).

Write $\pi:\R^d\to\R^m$ for projection onto the first $m$ coordinates, and set $V=F^{-1}(E)=\pi(E)$.
Choose an open set $W\subset\R^m$ containing $V$ with $\Omega=\overline W$ compact in $U$.
Apply Lemma~\ref{lem:taylor} with this $\Omega$ and the fixed $C_0>3\sqrt m/4$ to obtain the constants $\delta_0$ and $C_1$ for the adapted boxes.

Fix a sufficiently small $0<\delta\le\delta_0$.
To construct the cutoffs in parameter space, choose a nonnegative $\varphi\in C_c^\infty((-3/4,3/4)^m)$ such that $\sum_{v\in\Z^m}\varphi(u-v)=1$ for every $u\in\R^m$.
Such a function is obtained by taking a nonnegative smooth bump equal to one on $[-1/2,1/2]^m$ and supported in $(-3/4,3/4)^m$, then dividing it by the positive, periodic sum of its integer translates.
For each $v\in\Z^m$, set $I_v=\delta(v+[-1/2,1/2]^m)$, with center $u_{I_v}=\delta v$, and define $\psi_{I_v}(u)=\varphi(u/\delta-v)$ for $u\in\R^m$.
Let $V_\delta=\{u\in\R^m:\operatorname{dist}(u,V)<\delta\}$, and let $\mathcal I_\delta$ consist of the cubes $I_v$ whose cutoff supports meet $V_\delta$.
Since every term that can be nonzero on $V_\delta$ is retained, $\sum_{I\in\mathcal I_\delta}\psi_I=1$ there.
Moreover, for every $I\in\mathcal I_\delta$ and $\nu\in\N_0^m$,
\begin{equation}\label{eq:fixed-scale-cutoff-properties}
 0\le\psi_I,\qquad \supp\psi_I\subset(3/2)I\subset B(u_I,C_0\delta),\qquad
 \norm{\partial^\nu\psi_I}_\infty\le\norm{\partial^\nu\varphi}_\infty\delta^{-|\nu|}.
\end{equation}
The derivative bound in \eqref{eq:fixed-scale-cutoff-properties} gives \eqref{eq:psi-bounds}.
The support inclusion in \eqref{eq:fixed-scale-cutoff-properties} also gives bounded overlap, since the grid centers have spacing $\delta$.
Also, each retained center is within $(1+3\sqrt m/4)\delta$ of $V$.
Thus, for sufficiently small $\delta$, all retained centers lie in $W\subset\Omega$ and $V_\delta\subset W$.

We next bound the number of retained cubes, since this is the quantity that enters the kernel estimates.
Let $\xi_1,\ldots,\xi_{\cov(E,\delta)}\in\R^d$ be the centers of a cover of $E$ by balls of radius $\delta$.
Because $\pi$ is $1$-Lipschitz, the balls in $\R^m$ of radius $\delta$ centered at $\pi(\xi_j)$ cover $V$.
Each retained center is therefore within $(2+3\sqrt m/4)\delta$ of some $\pi(\xi_j)$.
A ball of this radius contains at most a constant number of points of $\delta\Z^m$, depending only on $m$.
Consequently, with $n_\delta=|\mathcal I_\delta|$, we have
\begin{equation}\label{eq:number-grid}
 n_\delta\lesssim\cov(E,\delta).
\end{equation}
We now pass to frequency space.
For each $I\in\mathcal I_\delta$, let $Q_I=Q(u_I,\delta;C_1)\subset\R^d$ be the adapted box centered at $F(u_I)$.
Recall from \eqref{eq:bI} that its cutoff is
\[
 b_{Q_I}(\xi)=\psi_I(\pi(\xi))\,\rho\bigl(D_\delta^{-1}A_{u_I}^T(\xi-F(u_I))\bigr),\qquad \xi\in\R^d.
\]
Set $\chi_\delta=\sum_{I\in\mathcal I_\delta}b_{Q_I}\in C_c^\infty(\R^d)$.
By \eqref{eq:fixed-scale-cutoff-properties}, the parameter cutoffs satisfy the support and derivative hypotheses of Lemma~\ref{lem:cutoff-neighborhood}.
Since they sum to one on $V_\delta\subset U$, the lemma gives
\[
 \chi_\delta(\xi)=\sum_{I\in\mathcal I_\delta}b_{Q_I}(\xi)=1
 \qquad\text{for every $\xi$ in an open neighborhood of $E=F(V)$.}
\]
Set $K_\delta=\check\chi_\delta\in\Sch(\R^d)$.
All cubes have side length $\delta$, and the parameter supports have bounded overlap.
Thus we may use the Lorentz bound \eqref{eq:assembly-lorentz} and the $L^\infty$ bound \eqref{eq:assembly-basic} from Lemma~\ref{lem:assembly}.
Substituting the count \eqref{eq:number-grid} into these bounds gives
\begin{align}
 \norm{K_\delta}_{p',1}
 &\lesssim \cov(E,\delta)^{1-1/p}\delta^{P_{\mathcal M}/p}
 =A_\delta^{1-1/p}\delta^{[P_{\mathcal M}-\alpha(p-1)]/p},\label{eq:fixed-main}\\
 \norm{K_\delta}_\infty
 &\lesssim\cov(E,\delta)\delta^{P_{\mathcal M}}
 =A_\delta\delta^{P_{\mathcal M}-\alpha}.\notag
\end{align}
Since $p<p_{\mathcal M,\alpha}$, the power of $\delta$ in \eqref{eq:fixed-main} is positive.
Therefore, along the chosen sequence $\delta_n\downarrow0$, boundedness of $A_{\delta_n}$ gives $\norm{K_{\delta_n}}_{p',1}\to0$.
Since $\chi_{\delta_n}=1$ on a neighborhood of $E$, Lemma~\ref{lem:cutoff}(i) gives $f=0$ almost everywhere.
This proves part~(i).

\smallskip\noindent\textit{Part (ii).} Assume $p=p_{\mathcal M,\alpha}$ and recall the endpoint hypothesis \eqref{eq:zero-content}:
\begin{equation}\label{eq:covering-proof-zero-content}
 \liminf_{\delta\downarrow0}\delta^\alpha\cov(E,\delta)=0.
\end{equation}
By \eqref{eq:covering-proof-zero-content}, choose $\delta_n\downarrow0$ with $A_{\delta_n}\to0$.
Since $P_{\mathcal M}>m\ge\alpha$, we have $p>2$ and may use the cutoff construction and \eqref{eq:fixed-main} with this value of $p$.
At the endpoint the power of $\delta$ in \eqref{eq:fixed-main} is zero; nevertheless, $A_{\delta_n}\to0$ gives $\norm{K_{\delta_n}}_{p',1}\to0$.
As in part~(i), the cutoffs equal one near $E$, so Lemma~\ref{lem:cutoff}(i) gives $f=0$ almost everywhere.
This completes the proof of Theorem~\ref{thm:covering}.
\end{proof}

\section{Applications}\label{sec:applications}

In this section, we give three applications of our main results.
First, we obtain necessary conditions for Fourier integrability of measures supported on finite-type manifolds.
The same argument rules out extension estimates in a range of exponents.
By Plancherel, the uniqueness results also give obstructions to $L^2$ densities for convolution powers.

\subsection{Necessary conditions for Fourier integrability of measures}\label{sec:fourier-integrability}

Zhang \cite[\S6.4.2, Question~6.3]{Zhang26} asks which support conditions on a manifold guarantee Fourier integrability at a common finite exponent for all probability measures satisfying those conditions.
Our uncertainty principle gives a necessary condition.
Let $\mu$ be a nonzero finite measure with compact support $E\subset\mathcal M$, and assume the graph charts covering $E$ have the common rank profile in Theorem~\ref{thm:covering}.
If $\widehat\mu\in L^{p,\infty}(\R^d)$ for some $1<p<\infty$ and $0<\dimH E\le m$, then Corollary~\ref{cor:uncertainty} gives
\[
 (p-1)\dimH E\ge P_{\mathcal M}.
\]
In strong $L^p$, we have a strict necessary condition under a finite-measure assumption.
Namely, if $\HH^\alpha(E)<\infty$ for some $0<\alpha\le m$ and $\widehat\mu\in L^p(\R^d)$ for $1\le p<\infty$, then
\[
 (p-1)\alpha>P_{\mathcal M}.
\]
Indeed, $f=\check\mu=\widehat\mu(-\,\cdot)$ is nonzero with Fourier transform $\mu$, so Theorem~\ref{thm:finite-hausdorff} rules out $p\le p_{\mathcal M,\alpha}$.
In weak $L^p$, finite lower $\alpha$-covering content still gives $p\ge p_{\mathcal M,\alpha}$ by Theorem~\ref{thm:covering}(i).
If that content is zero, part~(ii) rules out equality.

This condition concerns global integrability at a fixed exponent.
In contrast, Algom and Khalil \cite[Theorem~1.1]{AlgomKhalil25} and Algom and Orponen \cite[Theorem~1.2]{AlgomOrponen26} show that for every $\varepsilon>0$ there is an exponent $p$ such that
\[
 \int_{B(0,R)} |\widehat\mu(\xi)|^p\,d\xi\lesssim_\varepsilon R^\varepsilon,\qquad R\ge1,
\]
for the classes of measures they consider.
Since $p$ may depend on $\varepsilon$ and the bound grows with $R$, these estimates do not give the global integrability addressed here.

On the moment curve in $\R^d$, let $\alpha=\dimH E>0$.
Our necessary condition is $p\ge 1+d(d+1)/(2\alpha)$.
For $d\ge3$, this is the threshold $p_\alpha$ of Duan, Pham, and Ryou \cite[Theorems~1.1--1.2 and Proposition~3.1]{DPR25}: they construct measures with $\alpha$-dimensional support and global $L^p$ Fourier integrability for $p>p_\alpha$, and exclude strong $L^{p_\alpha}$ under their covering hypothesis.
For $d=2$, their threshold $4/\alpha$ is larger than $1+3/\alpha$ when $\alpha<1$.

\subsection{Obstructions to Fourier extension estimates}\label{sec:extension-obstructions}

Let $\mu$ be a finite positive Borel measure supported on a compact set $E\subset\mathcal M$.
Fix exponents $1\le p<\infty$ and $1\le q\le\infty$.
In the adjoint formulation of Fourier restriction, one asks whether the extension estimate
\[
 \|\widehat{h\mu}\|_{L^{p,\infty}(\R^d)}\lesssim\|h\|_{L^q(\mu)}
\]
holds for all $h\in L^q(\mu)$, with a constant independent of $h$.
Here $h\mu$ is the measure with density $h$ relative to $\mu$.
Since $\mu$ is finite and $q\ge1$, $h\in L^1(\mu)$ and $h\mu$ is a finite complex measure.
As in the moment-curve argument of Duan, Pham, and Ryou \cite[Proposition~1.1]{DPR25}, we show that in the ranges below the left-hand side is infinite for every nonzero $h$.

\begin{corollary}\label{cor:restriction-obstruction}
Let $\mathcal M$ and $E$ satisfy the hypotheses of Theorem~\ref{thm:hausdorff}, with $\alpha=\dimH E>0$.
Let $\mu$ be a finite positive Borel measure supported on $E$, let $1\le q\le\infty$, and let $h\in L^q(\mu)$ be nonzero.
Then
\[
 \widehat{h\mu}\notin L^{p,\infty}(\R^d)
 \qquad\text{for }1\le p<p_{\mathcal M,\alpha}.
\]
At $p=p_{\mathcal M,\alpha}$, the same conclusion holds if the lower $\alpha$-covering content of $E$ is zero.
If $\HH^\alpha(E)<\infty$, then $\widehat{h\mu}\notin L^{p_{\mathcal M,\alpha}}(\R^d)$ at the endpoint.
\end{corollary}

\begin{proof}
The measure $h\mu$ is nonzero and supported on $E$.
Thus $f=\check{(h\mu)}$ is nonzero and $\widehat f=h\mu$, while $f(x)=\widehat{h\mu}(-x)$ for $x\in\R^d$.
Below the endpoint, Theorem~\ref{thm:hausdorff} rules out the stated weak integrability.
If the lower content is zero, Theorem~\ref{thm:covering}(ii) does so at the endpoint as well.
If $\HH^\alpha(E)<\infty$, Theorem~\ref{thm:finite-hausdorff} applied to $f$ gives the strong endpoint exclusion.
This proves the corollary.
\end{proof}

\subsection{Distributions and convolution powers}
\label{sec:convolution-powers}

Theorem~\ref{thm:hausdorff} applies to distributions supported on $E$, not just to measures.
For example, let $T\in\Sch'(\R^d)$ be a distributional derivative of a measure supported on $E$.
If $E$ satisfies the hypotheses of that theorem and $\check T$ is represented by a weak-$L^p$ function for $1\le p<p_{\mathcal M,\alpha}$, where $\alpha=\dimH E$, then the theorem gives $\check T=0$ and hence $T=0$.

For positive measures, the Fourier-integrability obstruction above has a convolution form.
Let $\mu$ be a finite positive Borel measure on $\R^d$, let $\ell\ge2$ be an integer, and write $\mu^{*\ell}$ for its $\ell$-fold convolution.
Since $\widehat{\mu^{*\ell}}=(\widehat\mu)^\ell$, Plancherel gives
\begin{equation}\label{eq:convolution-plancherel}
 \mu^{*\ell}\text{ has an }L^2(\R^d)\text{ density}
 \quad\Longleftrightarrow\quad
 \widehat\mu\in L^{2\ell}(\R^d).
\end{equation}
In the reverse direction, the inverse Fourier transform of $(\widehat\mu)^\ell$ is an $L^2$ function representing $\mu^{*\ell}$ as a distribution.
Thus we can rule out an $L^2$ density at a prescribed order $\ell$ by applying uniqueness at exponent $2\ell$.

\begin{corollary}[Obstruction to convolution regularization]
\label{cor:convolution}
Let $d,m,k$ be integers with $1\le m<d$ and $k\ge2$.
Let $\mathcal M\subset\R^d$ be an $m$-dimensional $C^k$ submanifold with constant rank profile $(r_1,\ldots,r_k)$ and $r_{k-1}<r_k=d$, as in Theorem~\ref{thm:hausdorff}.
Let $\mu$ be a nonzero finite positive measure supported on a compact set $E\subset\mathcal M$, and let $\ell\ge2$ be an integer.
Assume one of the following:
\begin{enumerate}
\item[(i)] $0<\alpha=\dimH E\le m$ and
\[
 \ell<\frac12\left(1+\frac{P_{\mathcal M}}{\alpha}\right);
\]
\item[(ii)] for some $0<\alpha\le m$,
\[
\HH^\alpha(E)<\infty,
 \qquad \ell\le\frac12\left(1+\frac{P_{\mathcal M}}{\alpha}\right).
\]
\end{enumerate}
Then $\mu^{*\ell}$ has no $L^2(\R^d)$ density.
In particular, if $\mu$ is a probability measure and $X_1,\ldots,X_\ell$ are independent with law $\mu$, then $X_1+\cdots+X_\ell$ has no $L^2$ density.
\end{corollary}

\begin{proof}[Proof of Corollary~\ref{cor:convolution}]
Let $\mu$, $E$, and the integer $\ell\ge2$ be as in the statement.
Suppose for a contradiction that $\mu^{*\ell}$ has an $L^2(\R^d)$ density.
Write $f=\check\mu$ on $\R^d$.
By \eqref{eq:convolution-plancherel} and reflection, $f\in L^{2\ell}(\R^d)$; moreover, $\supp\widehat f\subset E$ and $f(0)=\mu(\R^d)>0$.
To apply \eqref{eq:finite-hausdorff-uniqueness} at exponent $2\ell$, set $s=P_{\mathcal M}/(2\ell-1)\in(0,P_{\mathcal M})$.

In case~(i), $s>\dimH E$, so $\HH^s(E)=0$.

In case~(ii), $\HH^\alpha(E)<\infty$ by assumption and $s\ge\alpha$.
Hence $\HH^s(E)<\infty$ also in this case.
Since $2\ell=1+P_{\mathcal M}/s$, \eqref{eq:finite-hausdorff-uniqueness} gives $f=0$ almost everywhere.
Since $f$ is continuous, this contradicts $f(0)>0$. \par
This completes the proof of Corollary~\ref{cor:convolution}.
\end{proof}

Convolution powers also occur in sharp Fourier restriction: Oliveira e Silva and Quilodr\'an \cite[Theorem~1.1]{OQ20} compare them for projection measures on the paraboloid and its convex perturbations, while Carneiro, Negro, and Oliveira e Silva \cite[Eq.~(6.15)]{CNO24} use the two-fold convolution of surface measure on the sphere.
For Frostman measures on the parabola, Orponen, Puliatti, and Py\"or\"al\"a \cite[Theorem~1.1 and Remark~1.4]{OPP25} obtain finite Riesz energy for sufficiently high convolution powers.
In contrast, the global Fourier integrability of the measures constructed by Duan, Pham, and Ryou \cite[Proposition~3.1]{DPR25} gives actual $L^2$ densities at sufficiently high orders by \eqref{eq:convolution-plancherel}.

\section{Examples}\label{sec:examples}

In the first example, we show that the upper bound on $p$ in Theorem~\ref{thm:hausdorff} cannot be increased on the product of two cubic moment curves.
For each $0<\alpha\le2$ and $p>1+12/\alpha$, we construct a compact set of Hausdorff dimension $\alpha$ on this product that supports the Fourier transform of a nonzero weak-$L^p$ function.

In the second example, we compare Theorem~\ref{thm:hausdorff} with earlier uniqueness results on a full two-dimensional surface patch in $\R^{14}$.
On this patch, the rank profile gives the range $1\le p<22$, which is larger than the ranges obtained by applying the bounds discussed below.

\subsection{Sharpness on a product of cubic moment curves}\label{sec:product-sharpness}

We use the graph parametrization $F_B:\R^2\to\R^6$ given by
\[
 F_B(u,v)=(u,v,u^2,v^2,u^3,v^3).
\]
Its rank profile is $(2,4,6)$, so its volume exponent is $12$.

The product of two cubic moment curves provides counterexamples above the uniqueness threshold for $F_B$.
We take two copies of the same Ahlfors regular measure on a cubic curve: the dimensions of their supports add, while factorization of the Fourier transform preserves the strong-$L^p$ integrability of the factor.

Let $0<\alpha\le d$.
A nonzero finite positive Borel measure $\mu$ with compact support $E\subset\R^d$ is $\alpha$-Ahlfors regular if $\mu(B(\xi,r))\approx r^\alpha$ for $\xi\in E$ and all sufficiently small $r>0$, uniformly in $\xi$.

If $F:X\to Y$ is a Borel map between Euclidean sets and $\mu$ is a finite Borel measure on $X$, its pushforward is defined by $F_\#\mu(A)=\mu(F^{-1}(A))$ for Borel sets $A\subset Y$.
If $\mu$ is a nonzero $s$-Ahlfors regular measure with compact support $E$, then $\dimH E=s$ and $\cov(E,\delta)\approx\delta^{-s}$ for all sufficiently small $\delta>0$; see \cite[Theorem~5.7 and its proof]{Mattila95}.

For finite Borel measures, the integral defining the Fourier transform uses $\dd\mu(x)$ in place of $f(x)\dd x$.
We use this convention also for the inverse transform.

\begin{proposition}\label{prop:product-sharpness}
Let $0<\alpha\le2$ and set $p_\alpha=1+12/\alpha$.
For every $p_\alpha<p<\infty$, there exist a compact set $E\subset F_B([0,1]^2)$ and a nonzero finite positive $\alpha$-Ahlfors regular measure $\mu$ with support $E$ such that
\[
 \dimH E=\alpha,\qquad
 \cov(E,\delta)\approx\delta^{-\alpha},\qquad
 0\ne\widecheck{\mu}\in L^{p,\infty}(\R^6).
\]
The covering estimate holds for all sufficiently small $\delta>0$.
The measure may depend on $p$.
\end{proposition}

\begin{proof}[Proof of Proposition~\ref{prop:product-sharpness}]
Fix $0<\alpha\le2$ and $p_\alpha<p<\infty$, and choose $p_\alpha<r<p$.
Let $\gamma:[0,1]\to\R^3$ be the cubic moment curve, $\gamma(t)=(t,t^2,t^3)$.

If $0<\alpha/2<1$, Duan, Pham, and Ryou \cite[Proposition~3.1 and proof of Theorem~1.2]{DPR25} give a nonzero finite positive $\alpha/2$-Ahlfors regular measure $\nu$ supported on $\gamma([0,1])$ with $\widehat\nu\in L^r(\R^3)$, since $r>1+12/\alpha$.
By reflection, $\widecheck{\nu}\in L^r(\R^3)$.

If $\alpha/2=1$, take $\nu=\gamma_\#(dt|_{[0,1]})$.
This measure is $1$-Ahlfors regular since $\gamma$ is bi-Lipschitz on $[0,1]$.
The convergence estimate of Arkhipov, Karatsuba, and Chubarikov \cite[Theorem~5]{AKC80} gives $\widecheck{\nu}\in L^r(\R^3)$ since $r>p_2=7$.

With $\nu$ chosen, let $\Pi:\R^6\to\R^6$ be the coordinate permutation that sends $(z_1,z_2,z_3,w_1,w_2,w_3)$ to $(z_1,w_1,z_2,w_2,z_3,w_3)$.
Set
\[
 \mu=\Pi_\#(\nu\otimes\nu),\qquad E=\supp\mu\subset F_B([0,1]^2).
\]
For $(z,w)\in\supp\nu\times\supp\nu$ and small $\rho>0$,
\[
 B(z,\rho/\sqrt2)\times B(w,\rho/\sqrt2)
 \subset B_{\R^6}((z,w),\rho)
 \subset B(z,\rho)\times B(w,\rho).
\]
The regularity of each factor and the isometry $\Pi$ therefore give $\mu(B(\xi,\rho))\approx\rho^\alpha$ for $\xi\in E$.
The preceding Ahlfors regularity bounds give $\dimH E=\alpha$ and $\cov(E,\delta)\approx\delta^{-\alpha}$.

Before the coordinate permutation, write $(x,y)\in\R^3\times\R^3$ for the physical variables.
The inverse transform factorizes:
\begin{equation*}
 \widecheck{(\nu\otimes\nu)}(x,y)=\widecheck{\nu}(x)\widecheck{\nu}(y),
 \qquad
 \norm{\widecheck{(\nu\otimes\nu)}}_r^r=\norm{\widecheck{\nu}}_r^{2r}<\infty.
\end{equation*}
Tonelli's theorem and invariance under $\Pi$ now give $f=\widecheck{\mu}\in L^r(\R^6)$ with $\widehat f=\mu$.
Since $\mu$ is finite, $f$ is bounded, and $f(0)=\nu(\R^3)^2>0$.
It follows that $0\ne f\in L^p(\R^6)\subset L^{p,\infty}(\R^6)$, as $p>r$.
This completes the proof of Proposition~\ref{prop:product-sharpness}.
\end{proof}

Since $P_{\mathcal M}=12$ on $F_B$, the examples rule out any increase of the covering threshold for $0<\alpha\le2$, even under an all-scale covering bound.
They also rule out any increase of the Hausdorff threshold throughout $0<\alpha\le2$.
Their lower $\alpha$-covering content is positive, so the endpoint criterion in Theorem~\ref{thm:covering}(ii) does not apply.

\subsection{A full surface with delayed rank gains}\label{sec:whole-surface-comparison}

We return to the surface from the introduction to see how its rank profile improves the uniqueness bound.
Let $I=[1/4,1/2]$ and write
\[
\begin{aligned}
 F_D(u,v)=\bigl(&u,u^2,u^3,u^4,u^6,\ v,v^2,v^3,v^4,v^6,\
                 uv,u^3v,uv^3,u^3v^3\bigr).
\end{aligned}
\]
Set $\mathcal M=F_D((0,1)^2)$ and $E=F_D(I^2)$.
The projection onto the $u$ and $v$ coordinates is the inverse of $F_D$ on $E$, and $F_D$ is Lipschitz on $I^2$.
Thus $E$ is bi-Lipschitz to $I^2$ and $\dimH E=2$.

\paragraph{Computing the rank profile.}
Fix $(u_0,v_0)\in(0,1)^2$ and write $x=u-u_0$ and $y=v-v_0$.
The rank $r_j$ is the rank of the matrix whose rows are the Taylor coefficients of the fourteen coordinate functions in degrees $1$ through $j$.
Indeed, each column is a partial derivative of $F_D$ divided by the corresponding multi-index factorial.
We compute this rank by changing the row basis after subtracting the constant terms.

In the pure $u$ block, the first four coordinates can be replaced by $x,x^2,x^3,x^4$.
Subtracting the Taylor polynomial of $u^6$ through order four leaves $6u_0x^5+x^6$, whose lowest-degree term is nonzero because $u_0>0$.
Similarly, the pure $v$ block gives $y,y^2,y^3,y^4$ and $6v_0y^5+y^6$.

For the mixed block, subtracting pure terms and constants from $uv,u^3v,uv^3,u^3v^3$ gives the four products of $x,(u_0+x)^3-u_0^3$ with $y,(v_0+y)^3-v_0^3$.
Subtracting the linear term from each cubic factor replaces these four products by
\[
 xy,\qquad (3u_0x^2+x^3)y,\qquad x(3v_0y^2+y^3),\qquad
 (3u_0x^2+x^3)(3v_0y^2+y^3).
\]
Their lowest-degree terms are nonzero multiples of $xy,x^2y,xy^2,x^2y^2$, respectively.
In particular, all fourteen lowest-degree monomials are distinct.
After truncation at degree $j$, the rows whose lowest degree exceeds $j$ vanish, and the others are independent: in any nontrivial linear combination, the terms of smallest degree cannot cancel.
Thus the new directions at each order are counted by
\[
\begin{array}{c|ccccc}
 j&1&2&3&4&5\\ \hline
 r_j-r_{j-1}&2&3&4&3&2\\
 r_j&2&5&9&12&14
\end{array}
\qquad (r_0=0).
\]
Since this computation holds at every base point, the surface has type five and the same rank profile throughout.
Consequently,
\[
 P_{\mathcal M}=\sum_{j=1}^5 j(r_j-r_{j-1})=2+6+12+12+10=42.
\]
As $\dimH E=2$, Theorem~\ref{thm:hausdorff} gives uniqueness in weak $L^p$ for $1\le p<22$.
Since $F_D$ is Lipschitz on $I^2$, $\HH^2(E)<\infty$.
Thus Theorem~\ref{thm:finite-hausdorff} also gives uniqueness in strong $L^p$ for $1\le p\le22$.

\paragraph{Comparison with earlier bounds.}
The dimension-based bound of Cao, Li, Miao, and Pang uses the sum of the degrees of the fourteen lowest-degree nonconstant monomials in two variables.
Those degrees sum to $2+6+12+20=40$, rather than $42$: in our example two of the five directions available at order four in the dimension count first appear at order five.
Deodhar's theorem applies with ambient dimension $14$ and surface dimension $2$.
Its auxiliary integer can be $2$ or $3$, yielding $(84-4)/5=16$ and $(84-10)/4=37/2$, respectively.

To apply Bi and Sun's result, set $(a,b,c,e)=(u,u^3,v,v^3)$.
Up to permutation, $F_D(u,v)$ lists these four coordinates and all ten quadratic monomials in them.
Thus $E$ lies in the quadratic Parsell--Vinogradov manifold with four parameters, to which their bound $4+\binom{6}{5}=10$ applies.
Table~\ref{tab:surface-comparison} summarizes the ranges on this patch.
The final two rows distinguish the strict weak-$L^p$ range from the strong-$L^p$ endpoint.

\begin{table}[htbp]
\centering
\begin{tabular}{@{}lcc@{}}
\toprule
Result & Function space & Uniqueness range\\
\midrule
Cao--Li--Miao--Pang \cite[Proposition~7.13]{CLMP25} & $L^p$ & $1\le p<21$\\
Deodhar \cite[Theorem~1.1]{Deo26} & $L^p$ & $2\le p\le 37/2$\\
Bi--Sun \cite[Theorem~1.1]{BiSun26} & $L^p$ & $1\le p<10$\\
Theorem~\ref{thm:hausdorff} & $L^{p,\infty}$ & $1\le p<22$\\
 Theorem~\ref{thm:finite-hausdorff} & $L^p$ & $1\le p\le22$\\
\bottomrule
\end{tabular}
\caption{Uniqueness ranges for functions with Fourier support in the patch $E\subset\R^{14}$.}
\label{tab:surface-comparison}
\end{table}
Since $E$ is compact, a weak-$L^p$ function with Fourier support in $E$ belongs to $L^q$ for every finite $q>p$, by Lemma~\ref{lem:lorentz}(iii) for $p>1$ and \eqref{eq:bandlimited-weak-l1} for $p=1$.
Thus the Cao--Li--Miao--Pang bound already gives weak-$L^p$ uniqueness for $1\le p<21$ by taking $p<q<21$.
The rank profile gives the additional range $21\le p<22$.
The argument gives strong uniqueness at $p=22$.

\section*{Use of AI tools}
During the preparation of this manuscript, the authors used ChatGPT (OpenAI) to assist with language editing, LaTeX preparation, checking mathematical arguments, exploring examples, and finding references.
The authors reviewed the resulting material, made revisions where necessary, and take full responsibility for the content of the manuscript.

\end{document}